\documentclass{article}
\usepackage[textwidth=6.5in,textheight=8.75in]{geometry}
\usepackage{graphicx} 

\usepackage{framed,multirow}

\usepackage{amsfonts,amssymb,booktabs}
\usepackage{latexsym}
\usepackage{hyperref}

\usepackage{url}
\usepackage{xcolor}
\definecolor{newcolor}{rgb}{.8,.349,.1}
\usepackage{authblk}
\usepackage{amsmath}

\usepackage{xcolor}

\usepackage{algorithmic}
\usepackage{algorithm,amsfonts,amssymb,amsthm}
\usepackage{mathrsfs}
\usepackage{appendix}
\usepackage{ascmac}

\usepackage{cleveref} 
\creflabelformat{equation}{#2(#1)#3}
\crefname{equation}{}{}

\crefname{figure}{Fig.}{}
\crefname{table}{Table}{}
\crefname{section}{Section}{}

\newtheorem{theorem}{Theorem}

\newtheorem{lemma}[theorem]{Lemma}

\newtheorem{proposition}[theorem]{Proposition}
\newtheorem{remark}[theorem]{Remark}

\def \v{\mathbf{v}}
\def \r{\mathbf{r}}

\def \m{\mathbf{m}}

\def \rp{\mathbf{r}_{p}}
\def \rq{\mathbf{r}_{q}}

\def \nq{\mathbf{n}_q}
\def \n{\mathbf{n}}

\def \K{K}
\def \V{V}

\def \rhohat{\hat{\boldsymbol{\rho}}}
\def \rhotil{\tilde{\boldsymbol{\rho}}}

\newcommand{\R}{{\mathbb{R}}}

\title{Layer potential quadrature on manifold boundary elements with constant densities for Laplace and Helmholtz kernels in $\R^3$}
\author[1]{Shoken Kaneko\footnote{kaneko60@umd.edu}}
\author[1]{Ramani Duraiswami\footnote{ramanid@umd.edu}}
\affil[1]{Department of Computer Science and Institute for Advanced Computer Studies, University of Maryland, College Park, MD 20742, USA.}

\begin{document}

\maketitle

\begin{abstract}
A method is proposed for evaluation of single and double layer potentials of the Laplace and Helmholtz equations on piecewise smooth manifold boundary elements with constant densities. 
The method is based on a novel two-term decomposition of the layer potentials, derived by means of differential geometry. 
The first term is an integral of a differential 2-form which can be reduced to contour integrals using Stokes’ theorem, while the second term is related to the element curvature. 
This decomposition reduces the degree of singularity and the curvature term can be further regularized by a polar coordinate transform. 
The method can handle singular and nearly singular integrals. 
Numerical results validating the accuracy of the method are presented for all combinations of single and double layer potentials, for the Laplace and Helmholtz kernels, and for singular and nearly singular integrals. 
\end{abstract}


\section{Introduction}
Boundary element methods (BEM) are widely used for solving partial differential equations arising in science and engineering. 
In the classical BEM, the boundary of the problem domain is typically represented using polygon meshes composed of piecewise flat boundary elements. 
This simple representation of the geometry allowed the development of efficient analytical methods tailored for flat boundary elements, e.g.,~\cite{GUMEROV2023recursive, kaneko2023recursive, lenoir2012evaluation, gumerov2023analytical}. 
Methods capable of solving problems with geometries represented by piecewise manifold surfaces which are not necessarily piecewise flat are receiving attention due to their ability to represent the geometry of a wide variety of problems accurately or exactly, and thereby eliminating a source of discretization error~\cite{beer2020isogeometric, greengard2021fast}. 

Practical BEM solvers are composed of multiple building blocks including iterative linear system solvers, fast matrix-vector product evaluation routines using fast multipole methods, etc. 
One of the essential computation routines in the BEM is the numerical evaluation of layer potentials integrals required for computing the near field interactions. 
This task is nontrivial because the integrands can be singular or nearly singular. 
Standard quadrature schemes which are effective for integrating polynomials of limited degrees, e.g. Gauss-Legendre quadrature, are known to produce inaccurate results when the evaluation point is close to the element. 
Many techniques have been developed over the years to accurately evaluate boundary integrals in such cases~\cite{adelman2016computation, guiggiani1992general, hackbusch1994numerical, hayami1988quadrature, johnston2007sinh, klockner2013quadrature, lenoir2012evaluation, montanelli2022computing, wala2019fast, wala2020optimization, zhu2022high}. 
Ref. \cite{montanelli2022computing} provides a recent extensive survey on this subject. 
The approaches developed include 
singularity cancellation using coordinate transforms~\cite{hayami1988quadrature, hackbusch1994numerical, johnston2007sinh}, 
singularity subtraction~\cite{guiggiani1992general}, 
continuation approach~\cite{rosen1995continuation}, 
dimension reduction~\cite{lenoir2012evaluation, zhu2022high, kaneko2023recursive, gumerov2023analytical}, 
adaptive subdivision~\cite{adelman2016computation}, and 
quadrature by expansion~\cite{klockner2013quadrature,wala2019fast,wala2020optimization}. 
The authors have recently proposed analytical methods based on dimensionality reduction for both collocation~\cite{kaneko2023recursive} and Galerkin BEM~\cite{gumerov2023analytical}, tailored for flat boundary elements. 
Zhu and Veerapaneni~\cite{zhu2022high} recently introduced a method for Laplace layer potentials on high-order curved elements using dimensionality reduction via Stokes' theorem and quaternion algebra. 
This method exploits the fact that an exact differential form is available for the Laplace double layer potential. 
While the application of this method to the evaluation of the Laplace single layer potential was discussed in passing in~\cite{zhu2022high}, numerical results were only presented for the double layer potential case.   
A summary of related quadrature methods for layer potentials or their multipole expansions is shown in \cref{table:relatedMethods}. 
\begin{table}[h]
\centering
\begin{tabular}{|c|c|c|}
\hline
 Distance $\to$ & Singular / nearly singular & Far-field expansions \\
 Element type $\downarrow$&    &   \\
\hline
\hline 
Flat, constant & Lenoir \& Salles~\cite{lenoir2012evaluation},   & Gumerov, Kaneko \&  \\
 &  Gumerov \& Duraiswami~\cite{gumerov2021analytical}& Duraiswami~\cite{GUMEROV2023recursive} \\
\hline
Flat, high order & Newman~\cite{newman1986distributions}, & Newman~\cite{newman1986distributions},\\
  & Kaneko, Gumerov \&  & Kaneko \& Duraiswami~\cite{kaneko2023efficient}\\
  & Duraiswami~\cite{kaneko2023recursive} & \\
\hline
Curved, constant & {\bf Present work} & \\
\hline
Curved, high order & Zhu \& Veerapaneni~\cite{zhu2022high}, & \\
& Klöckner et al.~\cite{klockner2013quadrature}, & \\
& Rosen \& Cormack~\cite{rosen1995continuation} & \\
\hline
\end{tabular}
\caption{Summary of related quadrature methods for layer potentials or its multipole expansions based on analytical or dimensionality-reduction based evaluation. Methods for high order elements can be applied to constant elements, and methods for curved elements can be applied to flat elements. 
A summary of methods based on other approaches e.g. singularity subtraction/cancellation can be found in~\cite{montanelli2022computing}.
}
\label{table:relatedMethods}
\end{table}

In this work, we focus on nearly singular and singular layer potential evaluation and propose a method which supports both Laplace and Helmholtz kernels for both single and double layer potentials on manifold boundary elements for the special case of constant densities. 
The method is based on a decomposition of the layer potentials into two terms. The first term is an integral of a differential 2-form, which can be evaluated via one-dimensional contour integrals after applying Stokes' theorem on manifolds, while the second term with reduced singularity which is related with the curvature of the element. 
The singularity in the second term can be further reduced by the classical technique of polar coordinate transform, used in e.g. \cite{hayami1988quadrature, guiggiani1992general, hackbusch1994numerical}. 
Layer potentials with higher order densities are important but require further development which may involve a redesign of the set of basis functions to obtain convenient exact differential forms, as shown for the Laplace double layer case~\cite{zhu2022high}. 
In this work we instead focus on the constant element case and provide formulations and numerical results for both Laplace and Helmholtz kernels for both the single and double layer potentials. 
The accuracy of the proposed method was confirmed via element-level tests and also using an example benchmark problem for which an analytical solution is available.

\section{Boundary element method and layer potentials}
The boundary element method is extensively used for numerical solution of partial differential equations, e.g. the Helmholtz equation and the Laplace equation, respectively given by
\begin{equation}\begin{aligned}
 - k^2 u(\r) - \nabla^2 u(\r) = f(\r), \quad  - \nabla^2 u(\r) = f(\r), \quad \r \in \Omega \subset \R^3,
\end{aligned}\label{eq:HelmholtzEqn} \end{equation}
with wavenumber $k$, field $u$ in domain $\Omega \subset \R^3$, and source $f$. 
The weak form of \cref{eq:HelmholtzEqn} can be written in terms of single- and double layer potentials $V$, $K$~\cite{sauter2010boundary}: 
\begin{equation}\begin{aligned}
&\{( c_p \gamma_{0,p} + \K\gamma_{0,q} - \V\gamma_{1,q}) u\}(\rp) = \{N_0f\}(\rp),\\
&\{\V\psi\}(\rp) \equiv \int_{\rq \in \Gamma} G(\mathbf{r}_p,\mathbf{r}_q) \psi(\rq) d\Gamma, \quad
\{ \K \phi\}(\rp) \equiv \int_{\rq \in \Gamma} \frac{\partial G(\mathbf{r}_p,\mathbf{r}_q)}{\partial \n_q} \phi(\rq) d\Gamma,\\
\end{aligned}\label{eq:BIE_BF_withRobinBC}
\end{equation}
with $c_p=1/2$ on a smooth boundary, $\gamma_{0}$ and $\gamma_{1}$ the boundary trace and normal derivative operators, and $N_0$ the Newton potential operator, defined as: 
\begin{equation}\begin{aligned}
&\{ \gamma_{0,q} u \} (\rq) \equiv \lim_{\hat{\r}_q \in \Omega \to \rq \in \Gamma} u(\hat{\r}_q), \quad 
\{ \gamma_{1,q} u \} (\rq) \equiv \nq \cdot \nabla_q u(\rq), \quad \rq \in \Gamma=\partial\Omega, \\
&\{N_0f\}(\rp) = \int_{\rq \in \Omega} G(\rp,\rq) f(\rq) d\Omega, \quad \rp \in \R^3,\\
\end{aligned} \end{equation}
where $G(\rp,\rq)$ is the respective Laplace or Helmholtz Green function:
\begin{equation}\begin{aligned}
G_\mathrm{L}(\rp,\rq) = \frac{1}{4\pi r}, \quad G_\mathrm{H}(\rp,\rq) = \frac{e^{ikr}}{4\pi r}, \quad  r \equiv |\mathbf{r}_q-\mathbf{r}_p|.
\end{aligned} \end{equation}

In the BEM the boundary $\Gamma$ is discretized into surface boundary elements which can be either flat or curved, and which may exactly discretize the original geometry when the closed-form representation of the geometry is available. This applies to e.g. surfaces generated using computer-aided design (CAD) software. The layer potential integrals over these elements are evaluated to form the linear system of equations. 
The densities $\psi$, $\phi$ are approximated via local, typically polynomial, functions (also called \emph{shape functions}) with unknown coefficients which must be determined.
In the present work we assume that the boundary $\Gamma = \partial \Omega$ is a union of boundary elements $\Gamma = \bigcup_i S_i$, where each $S_i$ is a smooth oriented Riemannian submanifolds with a boundary~\cite{lee2018introduction} and has constant density. Geometrical singularities e.g. wedges or corners need to be removed by subdividing the surface before applying the proposed method. 


\section{Differential geometry preliminaries}
\subsection{Curvature of regular surfaces}
\begin{figure}[htbp]
    \centering
    \includegraphics[width=6cm]{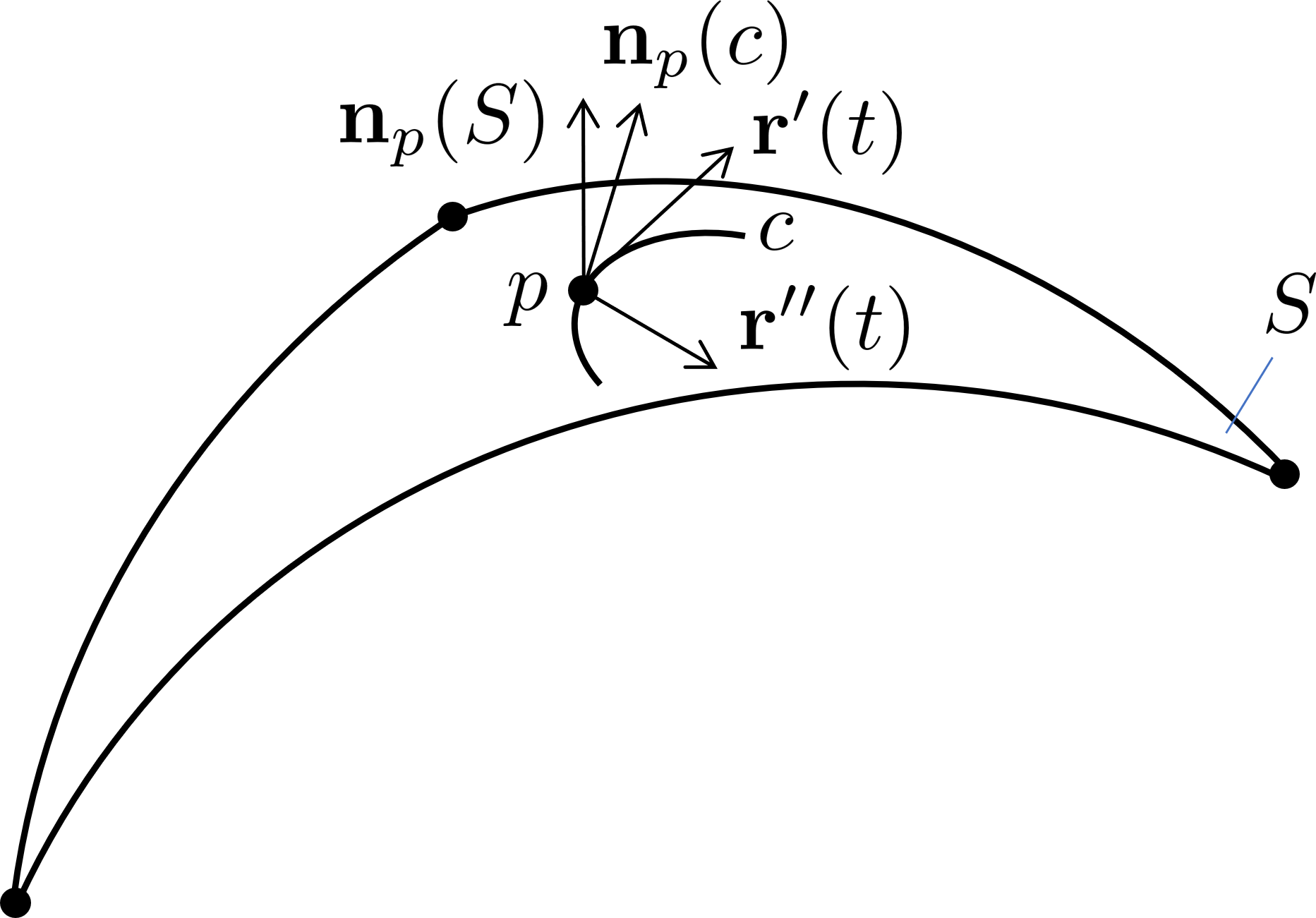}
    \includegraphics[width=6cm]{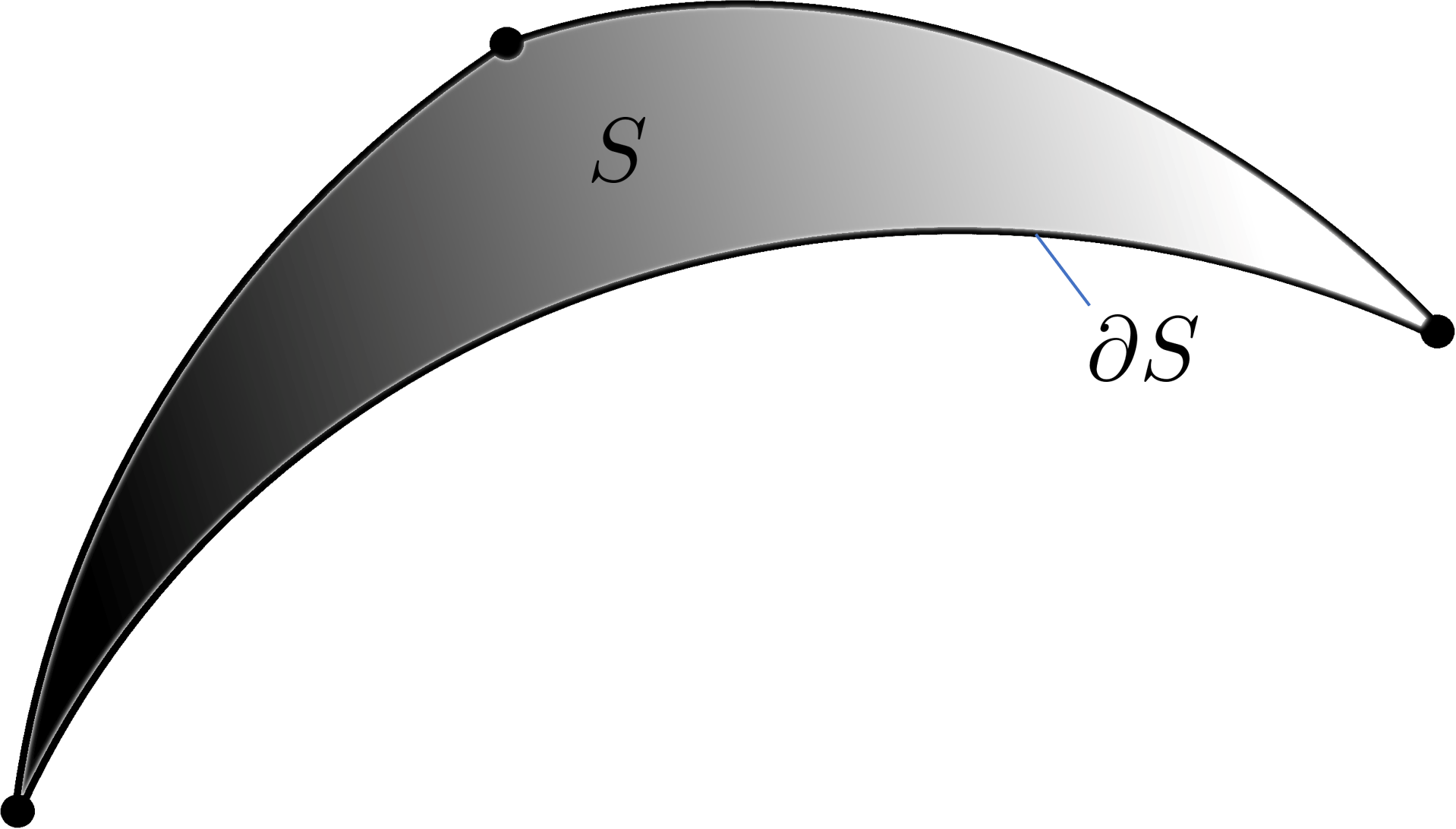}
    \caption{Left: The normal vectors of a regular surface $S$ and a regular curve $c$ at point $p$. Right: the surface $S$ and its contour $\partial S$.}
    \label{fig:element_surface_curve_normals}
    \end{figure}

In the differential geometry of curves and surfaces, various types of curvatures are defined. Here we briefly review the definition of the normal curvature, as it is central to the proposed method. The normal curvature $\kappa_{N}(p,c)$ of a regular curve $c$ on a regular surface $S$ at point $p \in c \subset S$ is defined as:
\begin{equation}\begin{aligned}
\kappa_{N}(p,c) &\equiv \kappa(p,c) \n_p(c) \cdot \mathbf{n}_p(S), \\
\kappa(p,c) &\equiv \frac{||\r''(t) \times \r'(t)||}{||\r'(t)||^3}, \quad \n_p(c)  \equiv \frac{\hat{\mathbf{t}}'(t)}{||\hat{\mathbf{t}}'(t)||}, \quad \mathbf{}
\end{aligned}\label{eq:curvatureOfCurve}
\end{equation}
 where $\kappa(p,c)$ is the curvature of curve $c$ at $p$, $p=\r(t) \in \R^3$ the parametrization of $c$, $\n_p(c)$ the unit normal vector of $c$ at $p$, $\mathbf{\hat{t}} = \r'(t)/||\r'(t)||$ the unit tangent vector of $c$ at $p$, and $\mathbf{n}_p(S)$ the unit normal vector of $S$ at $p$~\cite{do2016differential}. 
 This setup is illustrated in \cref{fig:element_surface_curve_normals} (left). 
 While $c$ can be any regular curve on $S$ passing through $p$, it is useful to consider the $c$ which is the cross section of $S$ with respect to a normal plane of $S$ at $p$ which is spanned by $\mathbf{n}_p(S)$ and a direction vector $\rhohat$, which is a tangent vector of $S$ at $p$. We refer to the normal curvature associated with direction $\rhohat$ as $\kappa_N(\rhohat)$. 
\subsection{Stokes' theorem on smooth manifolds}
 Another important tool available for the calculus on manifolds is Stokes' theorem~\cite{spivak2018calculus}:
\begin{equation}\begin{aligned}
\int_{S} \mathrm{d}\omega = \oint_{\partial S}\omega,
\end{aligned}\label{eq:stokesTheorem}\end{equation}
where $S$ is an oriented smooth submanifold with a boundary $\partial S$ as illustrated in \cref{fig:element_surface_curve_normals} (right), $\omega$ a differential form and $\mathrm{d}\omega$ the exterior derivative of $\omega$. 
Stokes' theorem allows the reduction of a surface integral into a contour integral, and can be used to evaluate integrals of exterior derivatives of exact differential forms. This applies to the Laplace double layer potential as demonstrated in~\cite{zhu2022high}. 

\section{Problem statement}\label{sec:problemStatement}
\begin{figure}[htbp]
    \centering
    \includegraphics[width=12cm]{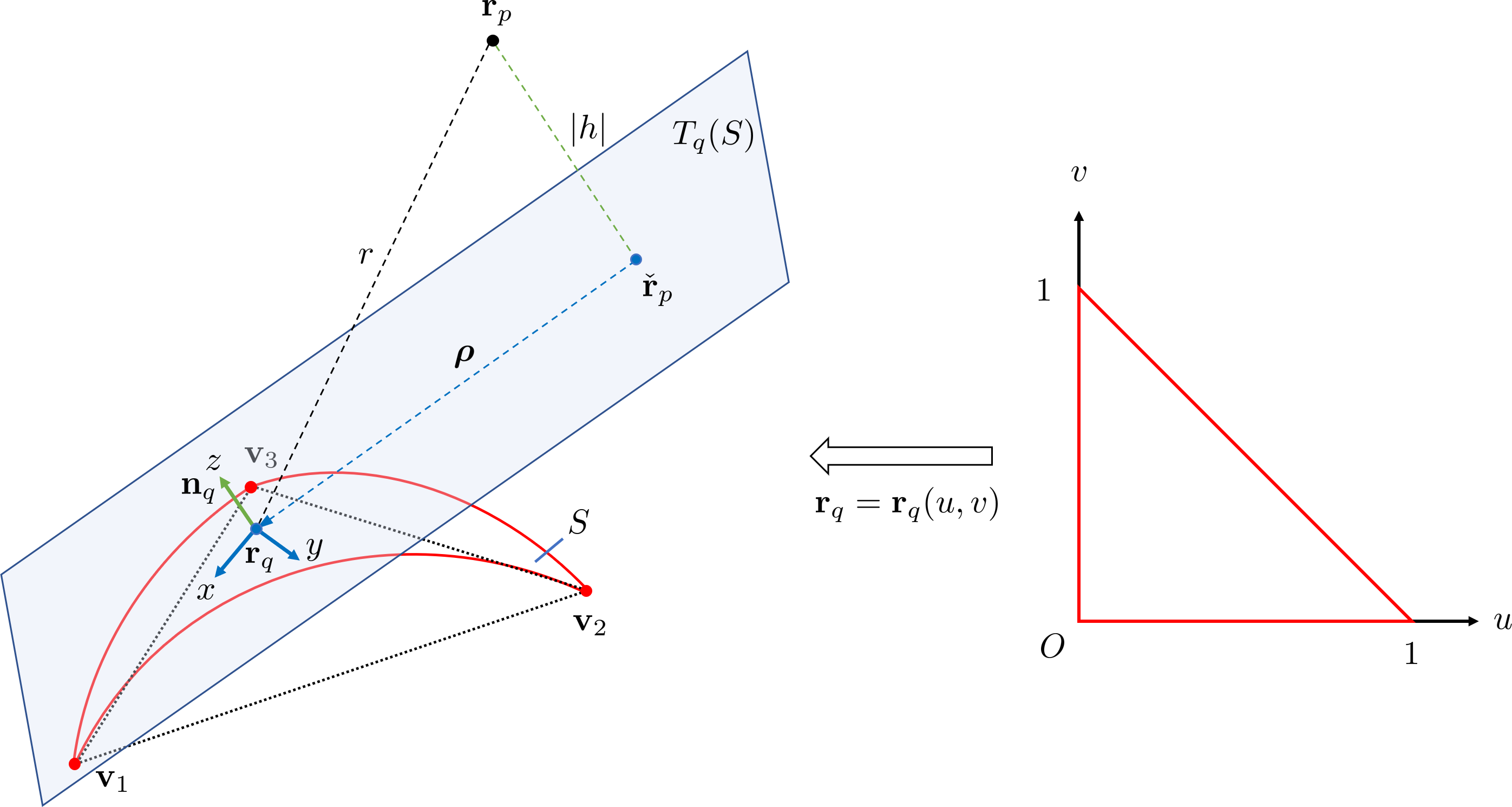}
    \caption{The manifold boundary element $S$ and the tangent plane $T_q(S)$ at $\rq$. Red, blue, and green lines indicate edges of the element, vectors parallel to the tangent plane $T_q(S)$, and vectors orthogonal to the tangent plane, respectively. $\nq$, $\rhohat \equiv \boldsymbol{\rho}/|\boldsymbol{\rho}|$ and $\rhotil \equiv \nq \times \rhohat$ can be used to construct a local orthogonal coordiante frame centered at $\rq$.}
    \label{fig:elementAndCoordinateFrames}
    \end{figure}

The setup under consideration is shown in \cref{fig:elementAndCoordinateFrames}. 
Let our boundary element $S$ be an oriented two-dimensional smooth Riemannian submanifold with a boundary~\cite{lee2018introduction} in $\R^3$, which is parametrized via the function $\rq(u,v) \in \R^3$ with variables $u, v \in \R$ defined on a reference triangle $\{u,v| 0\le u, 0\le v, u+v\le 1\}$. 
In the rest of the paper we will refer to this element simply as \emph{manifold element}. 
Let us denote the vertices of $S$ as $\v_1$, $\v_2$, and $\v_3$, 
and the unit normal vector and tangent plane at point $\rq$ as $\nq$ and $T_q(S)$, respectively.
Let us denote the normal plane at point $\rq$ spanned by $\nq$ and a tangent vector $\boldsymbol{t}$ as $N_q(\boldsymbol{t})$. 
The projection of an evaluation point $\rp$ onto a given tangent plane is denoted as $\check{\r}_p$. 
We also define $r \equiv |\mathbf{r}_q-\mathbf{r}_p|$, $\boldsymbol{\rho} \equiv \r_q-\check{\mathbf{r}}_p $, $\rho \equiv |\boldsymbol{\rho} |$, $h \equiv \nq \cdot (\rp-\rq )$. 
The goal is to develop a method to evaluate the single and double layer potentials over a given manifold element $S$ which is accurate in the nearly singular ($0 < \min_{\rq} r \ll 1$) and singular ($\min_{\rq} r=0$) cases. 

\section{Proposed method}

The following proposition provides a decomposition of the Green functions and their normal derivatives which allows the application of Stokes' theorem. It is a main result of our paper. 
\begin{proposition}[Decomposition of Green functions and their normal derivatives]\label{theorem:decomposition}
Green functions and their normal derivatives for the Laplace and Helmholtz equation can be decomposed into two terms as follows:
\begin{equation}\begin{aligned}
G_{\mathrm{K}}(\rp,\rq) &=  (\nabla_{\rq} \times  \mathbf{f}_{\mathrm{K}} )\cdot \nq  + \frac{1}{4\pi}\left(C_{\mathrm{K}}\kappa_N(\rhotil)+D_{\mathrm{K}}\kappa_N(\rhohat)\right),\\
\frac{\partial G_{\mathrm{K}}(\rp,\rq)}{\partial \nq} &=  (\nabla_{\rq} \times  \mathbf{f}_{\mathrm{K}}' )\cdot \nq + \frac{1}{4\pi}\left(C'_{\mathrm{K}}\kappa_N(\rhotil)+D'_{\mathrm{K}}\kappa_N(\rhohat)\right),\\
\end{aligned}\label{eq:greenFncDecomposition}\end{equation}
with $\mathrm{K}=\{\mathrm{L},\mathrm{H}\}$ the identifier whether the kernel is Laplace or Helmholtz, $\nq$ the unit normal vector at point $\rq$, $\kappa_N(\rhohat)$ and $\kappa_N(\rhotil)$ the normal curvature of the element at point $\rq$ for the normal planes spanned by $\nq$ and the tangent vectors $\rhohat \equiv \boldsymbol{\rho}/|\boldsymbol{\rho}|$ and $\rhotil \equiv \nq \times \rhohat$, respectively, $\mathbf{f}_{\mathrm{K}}$ and $\mathbf{f}'_{\mathrm{K}}$ the pseudo potential fields, $C_{\mathrm{K}}$, $C'_{\mathrm{K}}$, $D_{\mathrm{K}}$ and $D'_{\mathrm{K}}$ the weights of the curvatures defined as follows: 
\begin{equation}\begin{aligned}
&\mathbf{f}_{\mathrm{L}} \equiv \frac{\rho\rhotil}{4\pi(r+h)}, \quad \!
C_{\mathrm{L}} \equiv \frac{h}{r+h}, \quad \! 
D_{\mathrm{L}} \equiv \frac{r}{r+h}, \\ 
&\mathbf{f}'_{\mathrm{L}} \equiv \frac{\mathbf{f}_{\mathrm{L}}}{r}, \quad \!
C'_{\mathrm{L}} \equiv \frac{C_{\mathrm{L}}}{r}, \quad \!
D'_{\mathrm{L}} \equiv \frac{D_{\mathrm{L}}}{r}, \\
&\mathbf{f}_{\mathrm{H}} \equiv \frac{e^{ikr}-e^{ikh}}{4\pi ik \rho}\rhotil, \quad
C_{\mathrm{H}} \equiv \frac{h(e^{ikr}-e^{ikh})}{ik\rho^2}, \quad 
D_{\mathrm{H}} \equiv e^{ikh} - \frac{h(e^{ikr}-e^{ikh})}{ik\rho^2}, \\
&\mathbf{f}'_{\mathrm{H}} \equiv \frac{r e^{ikh} - h e^{ikr}}{4\pi r \rho}\rhotil, \quad \!\!\!\!
C'_{\mathrm{H}} \equiv \frac{h\left(r e^{i k h}-h e^{i k r}\right)}{r \rho^2}, \quad \!\!\!\!
D'_{\mathrm{H}} \equiv \frac{r e^{i k r}-h e^{i k h}}{ \rho^2}-i k e^{i k h}. \\
\end{aligned}\label{eq:f_C_D_expressions}\end{equation}
\end{proposition}
\begin{proof}
See \cref{appendix_proof_decomposition}.    
\end{proof} 
We refer to the first and second term in decomposition \cref{eq:greenFncDecomposition} as the \emph{Stokes term} and the \emph{Curvature term}, respectively. 
The Stokes term offers a differential 2-form: $\mathrm{d}\omega = (\nabla_{\rq} \times  \mathbf{f}_{\mathrm{K}} )\cdot \nq  dS$, whose integral can be reduced to a contour integral of a differential 1-form $\omega = \mathbf{f}_{\mathrm{K}}\cdot d\mathbf{l}$ due to Stokes' theorem with $d\mathbf{l}$ the line element vector along the contour $\partial S$.
Hence, the layer potentials now can be expressed as follows:
\begin{equation}\begin{aligned}
\{V_{\mathrm{K}}\}(\rp) &= \oint_{\partial S} \mathbf{f}_{\mathrm{K}} \cdot d\mathbf{l} + \frac{1}{4\pi} \int_S \left(C_{\mathrm{K}}\kappa_N(\rhotil)+D_{\mathrm{K}}\kappa_N(\rhohat)\right) dS,\\
\{K_{\mathrm{K}}\}(\rp) &= \oint_{\partial S} \mathbf{f}'_{\mathrm{K}} \cdot d\mathbf{l} + \frac{1}{4\pi} \int_S \left(C'_{\mathrm{K}}\kappa_N(\rhotil)+D'_{\mathrm{K}}\kappa_N(\rhohat)\right) dS,\\
\end{aligned}\label{eq:evaluateLayerPotentialUsingDecomposition}\end{equation}
with subscript $\mathrm{K}=\{\mathrm{L},\mathrm{H}\}$ indicating the type of the kernel. 

\begin{remark}
The double layer potential of the Laplace kernel can be expressed as an integral of an exact form only and~\cite{zhu2022high} utilized this fact. On the other hand, the single layer potential was approximated in~\cite{zhu2022high} by a double layer potential with a modified density function. The efficacy of this approach was not discussed in~\cite{zhu2022high} and is unclear. 
We present numerical results for both the present method based on decomposition \cref{eq:greenFncDecomposition} and the method presented in~\cite{zhu2022high} in \cref{sec:results} for a comparison. 
\end{remark}

\begin{remark}
    It is interesting that the proposed method, which is derived from the perspective of differential geometry, resembles a feature of the continuation approach~\cite{rosen1995continuation} which also results in a decomposition of the integrand into two parts where one of the terms absorbs the ``worst part" of the singularity. The decomposition in the continuation approach is based on Taylor series expansions and the geometric meanings of the decomposed terms are not clear. In contrast, decomposition \cref{eq:greenFncDecomposition} offers two terms with clear geometric meanings associated with differential forms and the curvature of the element.  
\end{remark}

\begin{remark}
If the element is flat, the curvature term vanishes and the layer potential evaluation \cref{eq:evaluateLayerPotentialUsingDecomposition} reduces to the evaluation of just contour integrals. 
Furthermore, this case can be evaluated analytically using the RIPE method~\cite{kaneko2023recursive}. 
\end{remark}

\begin{figure}[htbp]
    \centering
    \includegraphics[width=12cm]{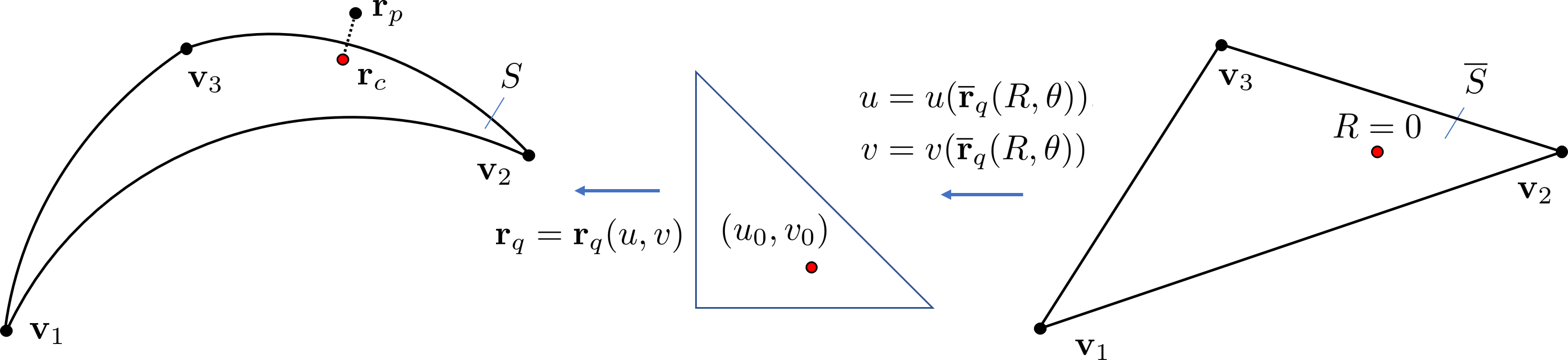}
    \caption{The coordinate transform used for the quadrature of the curvature term.}
    \label{fig:elementCoordinatesAndMaps}
    \end{figure}

As a consequence, the one-dimensional integral of the Stokes term over the curvilinear boundary can be evaluated using Gauss-Legendre quadrature, and methods for two-dimensional quadrature can be applied to the curvature term. 
For non-negative $h$, as $h\to0$ with $r\to0$, the curvature term has a regularity of $1/r^{n-1}$ as opposed to $1/r^n$ of the original integrand before the decomposition, where $n=1$ for the single layer potential and $n=2$ for the double layer potential. This means that the integral still contains a (near-) singularity in the double layer potential case. 
To further regularize the singularity, we employ the classical technique of polar coordinate transform around the singularity~\cite{hayami2005variable}, which is illustrated in  \cref{fig:elementCoordinatesAndMaps}. A polar coordinate system $(R,\theta)$ is set up on a flat surrogate element $\overline{S}$ whose vertices are identical to those of the original manifold element. 
Points on this surrogate element $\overline{S}$ are mapped to the reference triangle via an affine mapping, which are then mapped to the manifold element $S$ via the parametrization $\rq=\rq(u,v)$. 
Point $R=0$, i.e. the origin of the polar coordinates on $\overline{S}$, is chosen to be the point which maps onto $\r_c=\rq(u_0,v_0)$, the point on $S$ nearest to the evaluation point $\rp$. 
This mapping was also used in~\cite{hayami2005variable}, though~\cite{hayami2005variable} applies subsequent coordinate transformations. Here we only use the basic polar coordinate transform. 
\begin{remark}
    The Stokes term can be still nearly singular if $\rp$ is close to an edge of the element. This can be avoided by using nonconforming boundary elements where the collocation points are in the interior of the element. This is naturally satisfied in constant elements with center-panel collocation. 
\end{remark}
\begin{remark}
For negative values of $h$, the functions in \cref{eq:f_C_D_expressions} 
have a singularity at $h=-r$ which hinders the evaluation of the integrals using the presented approach. Such singularities arise when $\rp$ is in the inward normal bundle of the element. This can be resolved by utilizing respectively the symmetry and anti-symmetry of the layer potentials with respect to the exchange of variables $u$ and $v$ in the parametrization $\rq(u,v)$. For the single layer potential we have the symmetry:
\begin{equation}\begin{aligned}
\int_{u=0}^{1}\int_{v=0}^{1-u} G(\rp,\rq(u,v))J(u,v)dvdu = \int_{u=0}^{1}\int_{v=0}^{1-u} G(\rp,\rq(v,u))J(v,u)dvdu 
\end{aligned}\label{eq:symmetry_and_antisymmetry}\end{equation} 
with $J$ the Jacobian of the transform from $\rq$ to $(u,v)$ and for the double layer potential we have the anti-symmetry where this exchange results in a sign flip. 
This exchange of variables also flips the direction of the normal vectors and the sign of $h$.
 For a given evaluation point $\rp$, therefore, we can evaluate the same layer potentials using this symmetry property to avoid the singularity $h=-r$. 
 This may not be always possible, since the outward and inward normal bundles of the element can have a non-empty intersection and evaluation points in this intersection cannot avoid the singularity.  This can happen if the element is too curved. 
 In such cases, the element can be subdivided until the evaluation point $\rp$ can avoid the inward normal bundle in one of the parametrizations and the proposed method can be applied to the subdivided elements. 
 The geometry is illustrated in \cref{fig:normalBundles}. 
An example pseudo-code implementing the procedure to avoid singularities is listed in \cref{algo:proposedMethod}. 
Various optimization methods could be used for step 1 in \cref{algo:proposedMethod}. The Newton-Raphson method was used in our implementation. While this step adds additional computation cost, similar computation is needed in other methods for nearly singular integrals (e.g.~\cite{hayami2005variable}), where the first step is to find the point on the element closest to $\rp$. As Newton's method converges quickly its computational overhead is limited and can be practically considered constant per element. 
    \begin{figure*}[htb!]
    \centering
    \includegraphics[height=5cm]{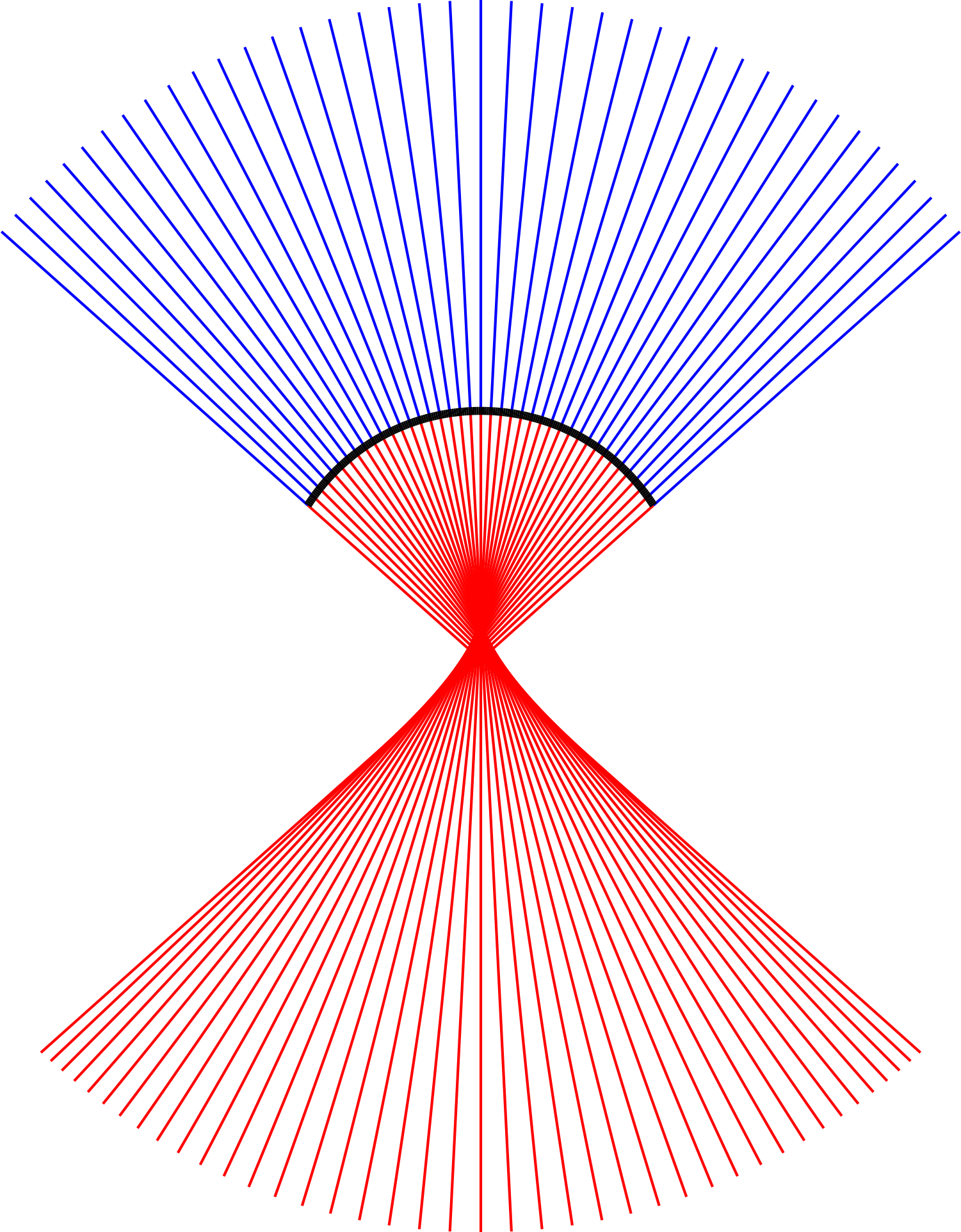}
    \includegraphics[height=5cm]{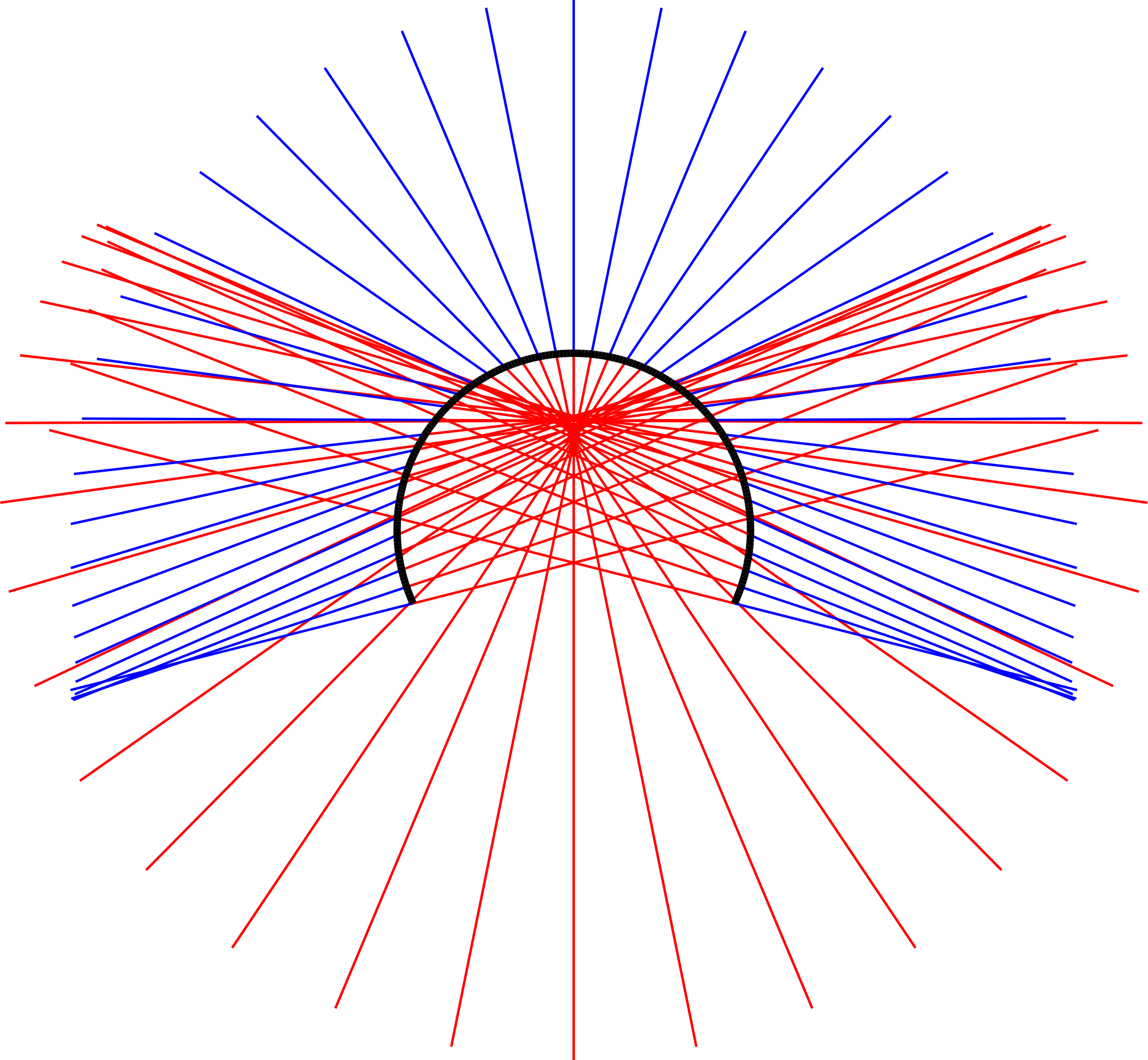}
    \caption{A 2D illustration of curved elements and their normal bundles. The element, their outward and inward normal bundles are drawn by black, blue, and red lines, respectively. Left: outward and inward normal bundles do not intersect. In this case, evaluation points in the inward normal bundle can avoid the singularity by flipping the parametrization. Right: outward and inward normal bundles intersect and if $\rp$ is in the intersection, the element has to be subdivided so that $\rp$ does not belong to an intersection of normal bundles.}
    \label{fig:normalBundles}
    \end{figure*}
\begin{algorithm}[htb!]
\caption{EvaluateLayerPotential($S$: element, $\rp$: evaluation point)
}\label{algo:proposedMethod}
\begin{algorithmic}
\STATE{1. Find $[\hat{h}_{\min},\hat{h}_{\max}]$, the range of $h/r$ over all points on the element.} 
\STATE{2. If $\hat{h}_{\min}=-1$ and $\hat{h}_{\max}=1$, subdivide the element such that $S = \bigcup_i s_i$. \\\quad \quad $I := \sum_i \mathrm{EvaluateLayerPotential}(s_i,\rp)$. Return $I$.}
\STATE{3. If $|\hat{h}_{\min} + 1| < |\hat{h}_{\max} - 1|$, flip the parametrization of $S$ from $\rq(u,v)$ to $\rq(v,u)$.} 
\STATE{4. Evaluate the integral $I$ using \cref{eq:evaluateLayerPotentialUsingDecomposition}. }
\STATE{5. If the target is the double layer potential and the parametrization \\\quad\quad was flipped in step 3., $I := -I$.}
\STATE{6. Return $I$.}
\end{algorithmic}
\end{algorithm}
\end{remark}

\section{Numerical evaluation}\label{sec:results}
\subsection{Element-level tests}
\subsubsection{Nearly singular case}\label{sec:nearlySingularResults}
\begin{figure*}[htbp]
    \centering
    \includegraphics[width=6cm]{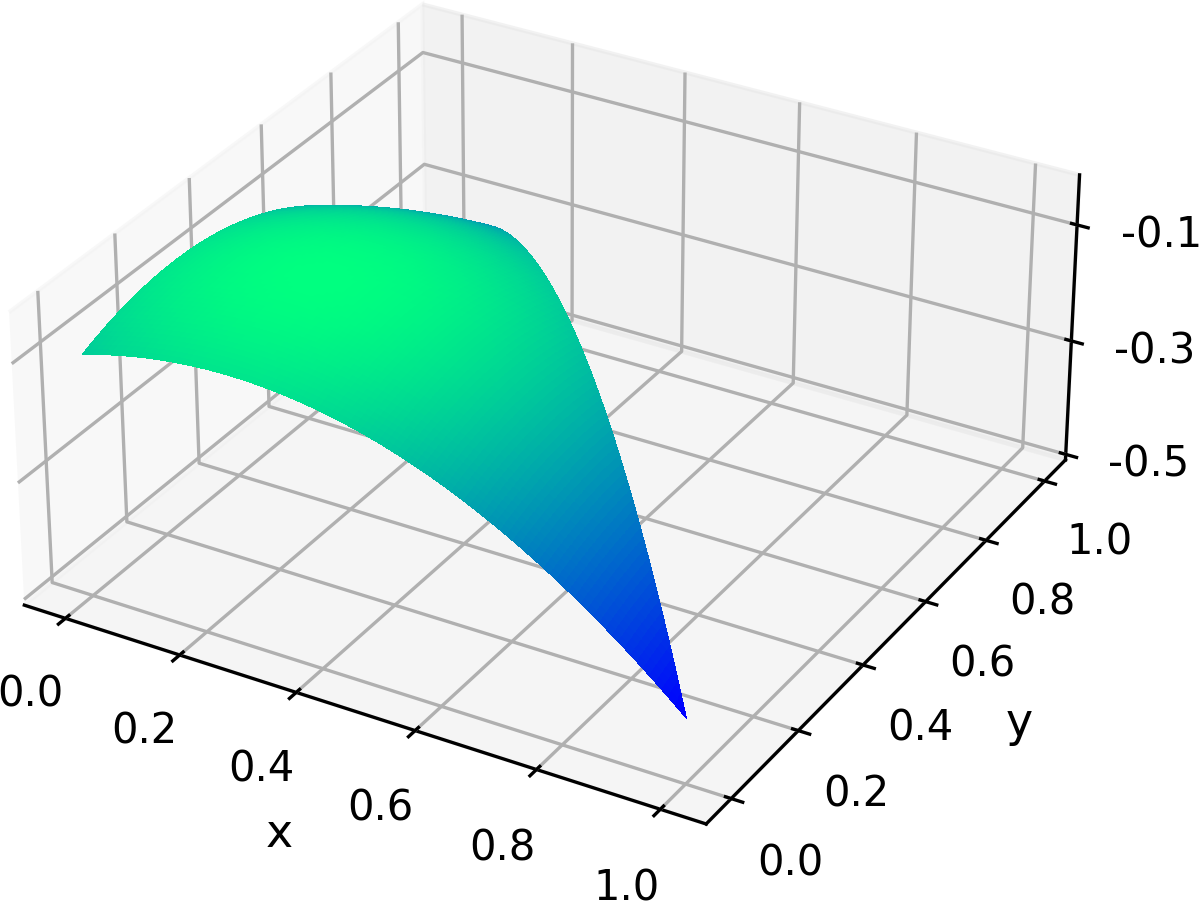}
    \includegraphics[width=6cm]{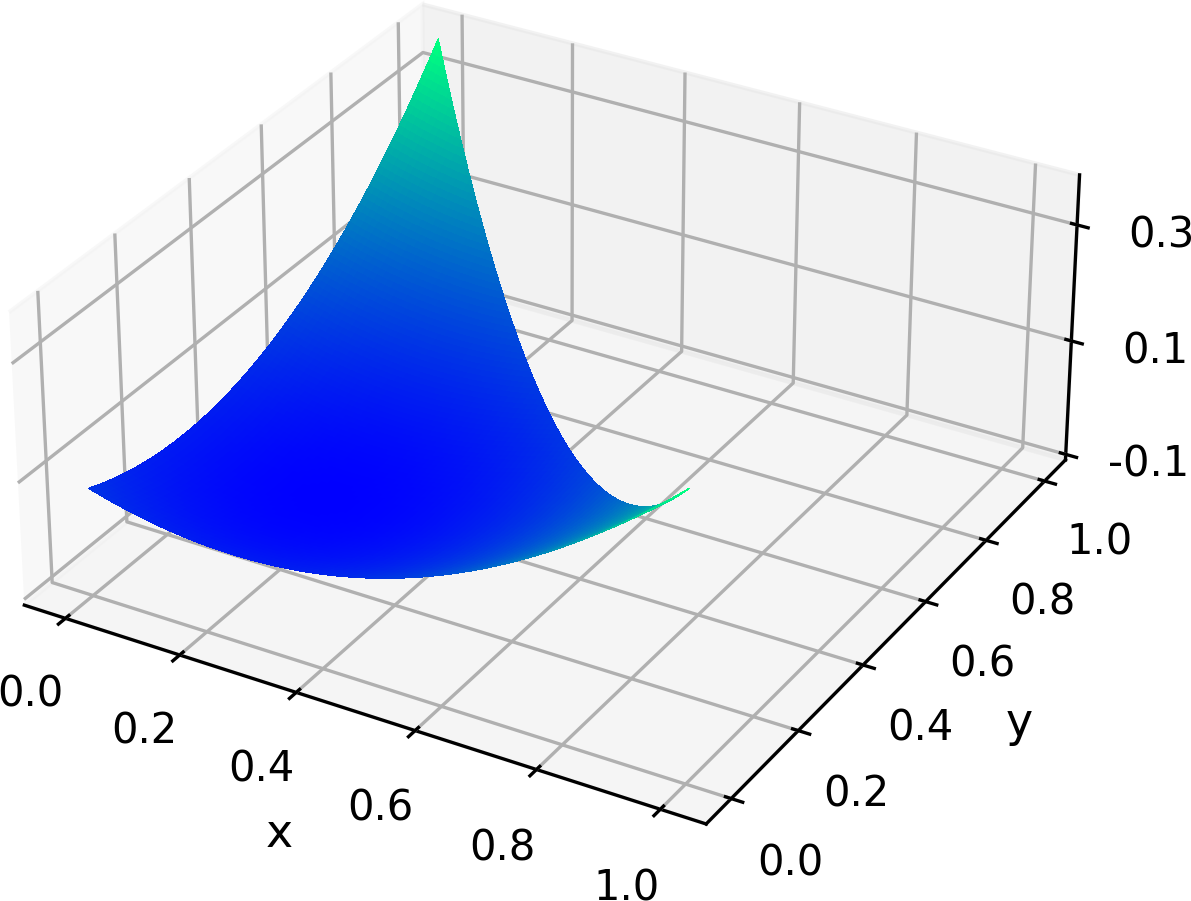}
    \caption{The manifold boundary elements with positive and negative curvature used in the numerical tests. \emph{Element 1} (left, $\sigma=-0.6$) and \emph{element 2} (right, $\sigma=0.6$). The colors are for visual aid.}
    \label{fig:testElements}
    \end{figure*}

\begin{figure*}[htbp]
    \centering
    \includegraphics[width=6cm]{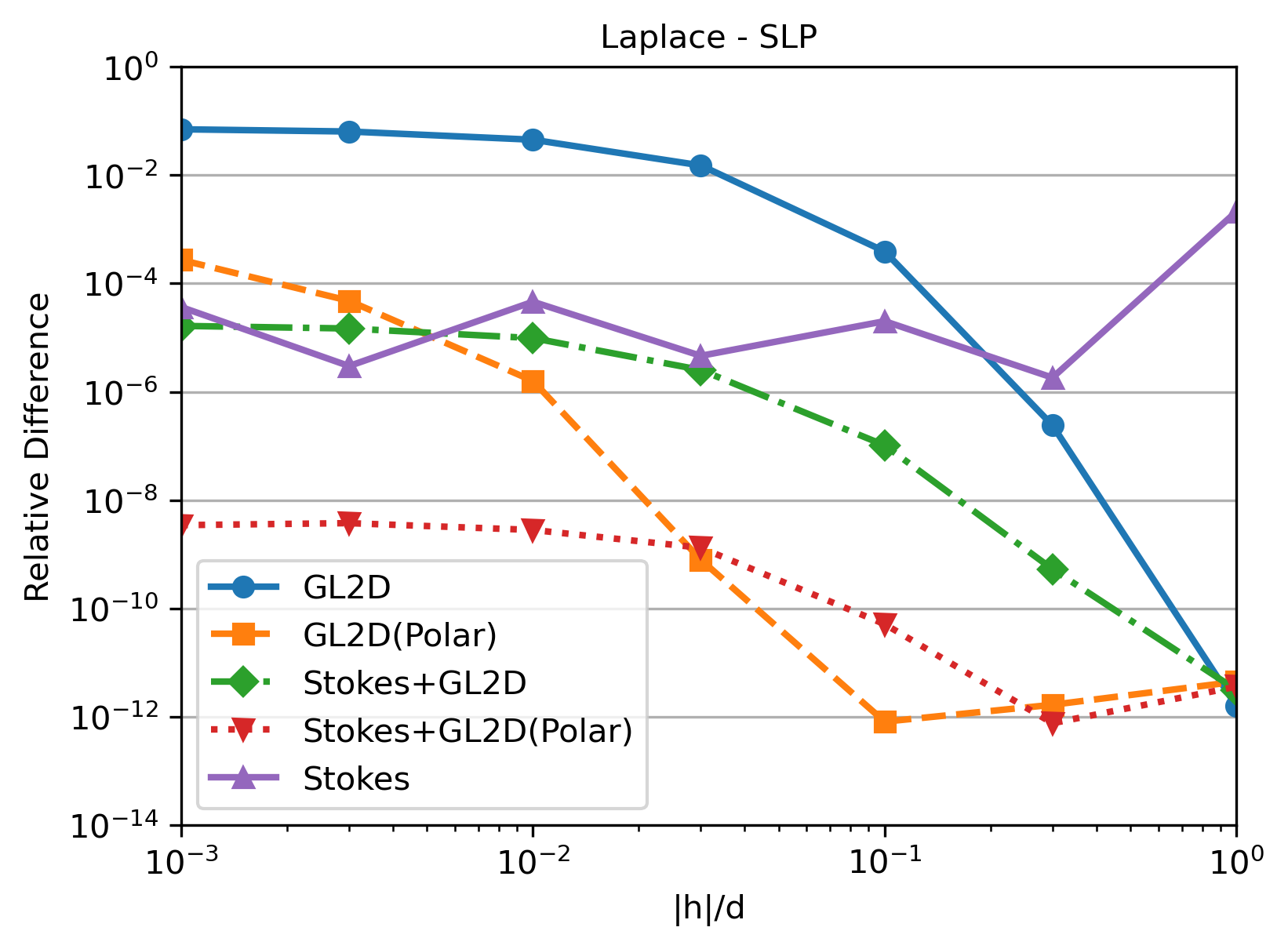}
    \includegraphics[width=6cm]{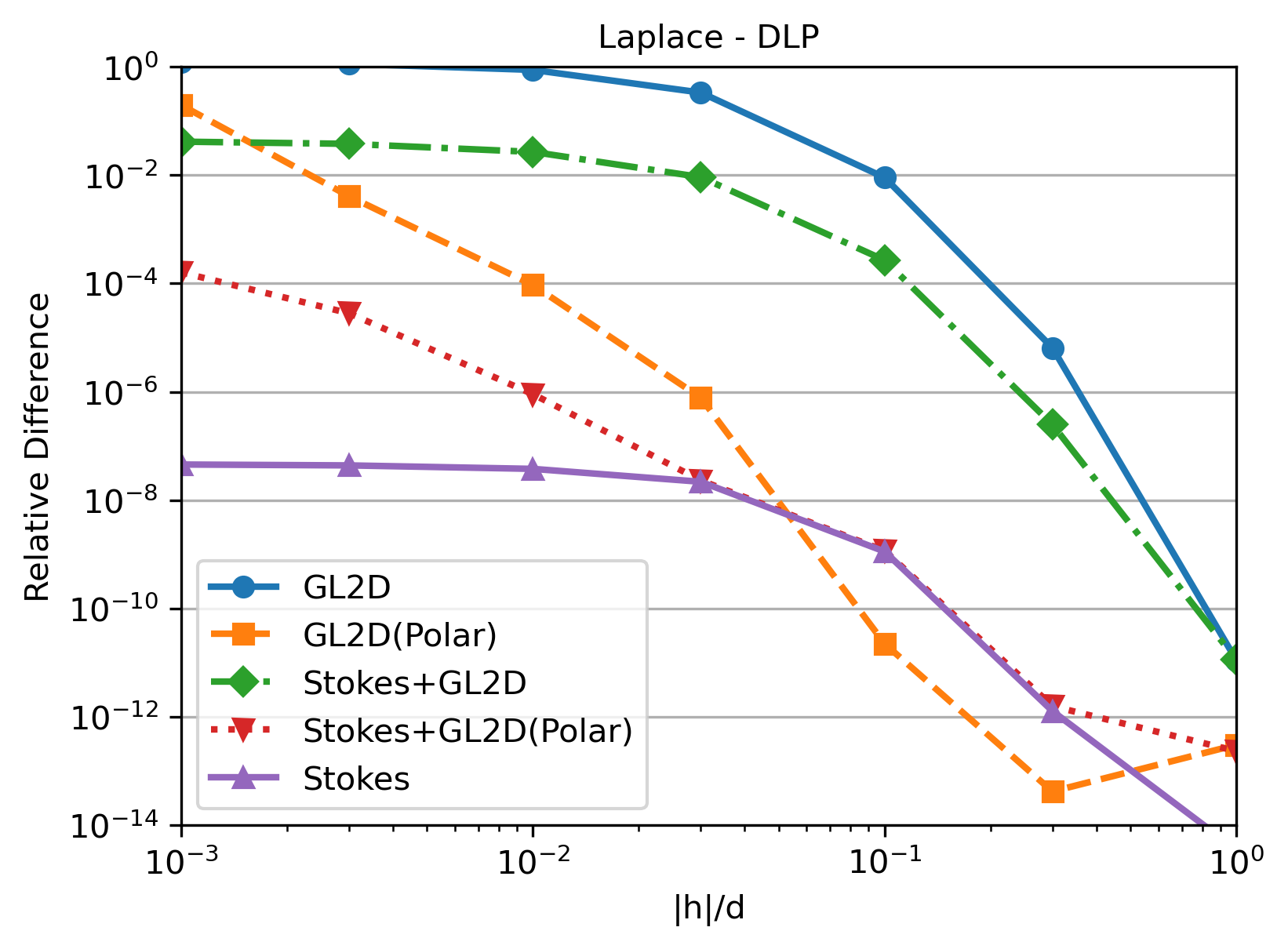}
    \includegraphics[width=6cm]{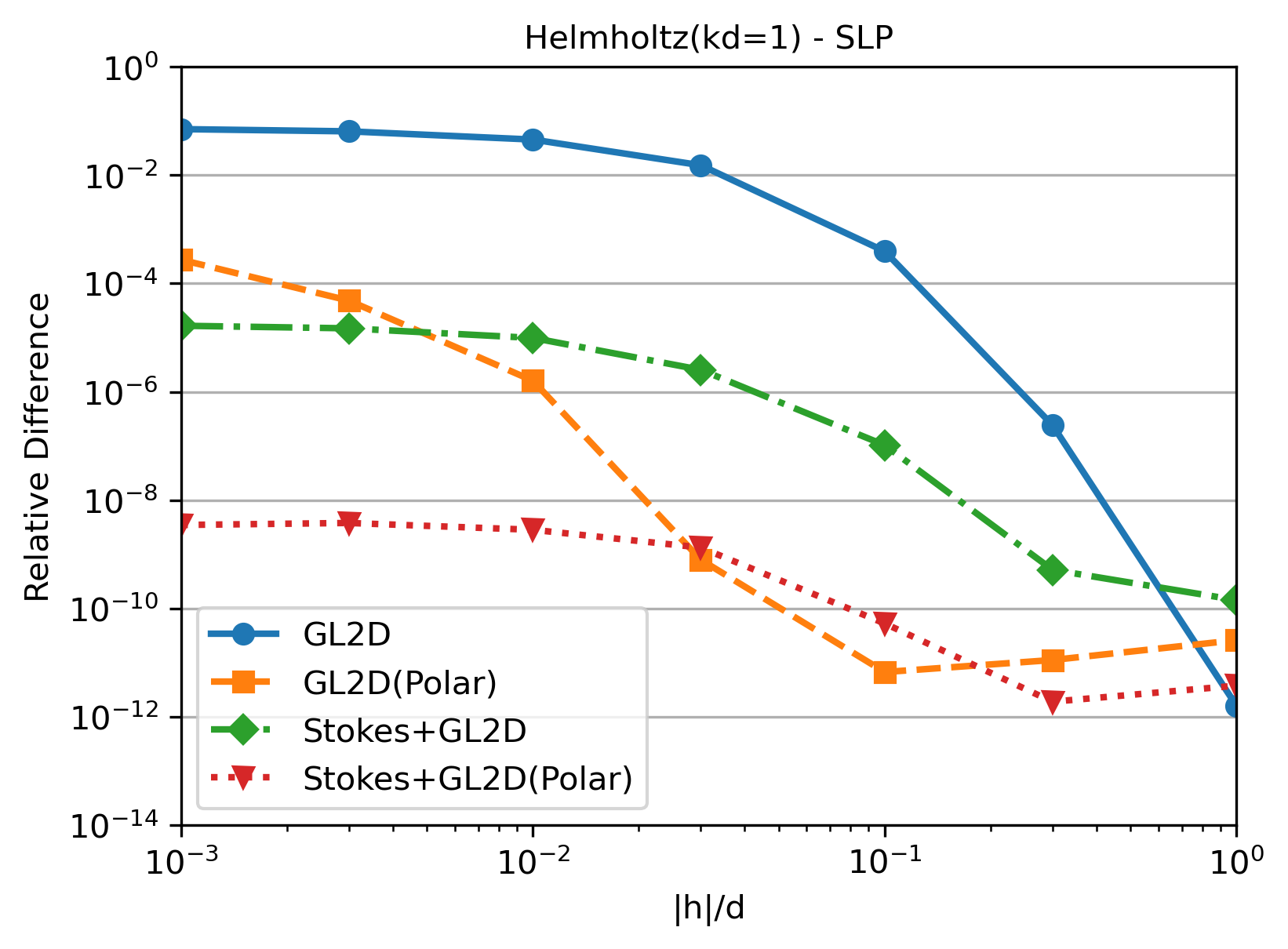}
    \includegraphics[width=6cm]{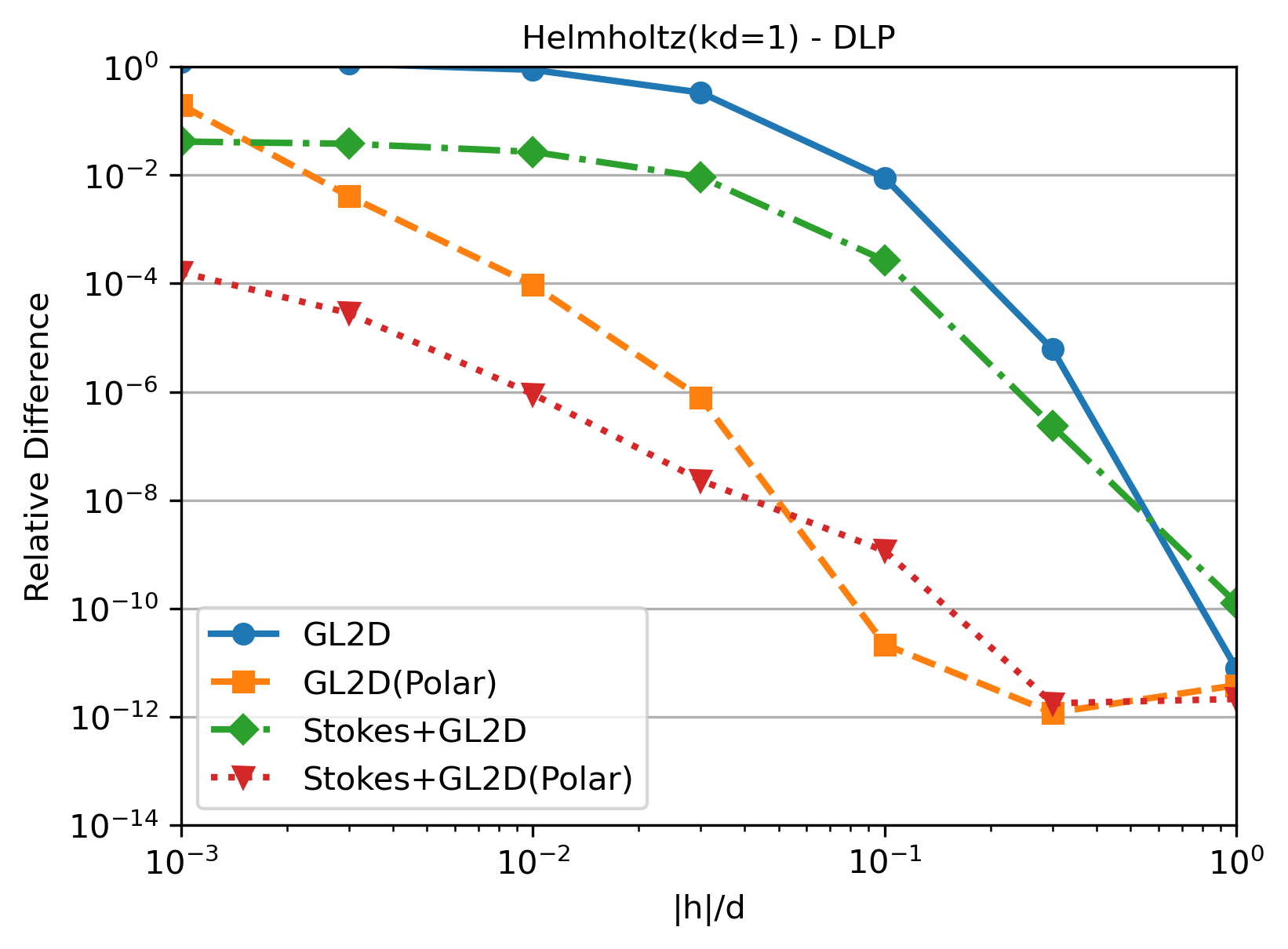}
    \caption{Relative difference of the layer potentials on \emph{element 1} for various methods and the reference adaptive Gauss-Kronrod quadrature. Results for single layer potential with Laplace kernel (top left), double layer potential with Laplace kernel (top right), single layer potential with Helmholtz kernel (bottom left), and double layer potential with Helmholtz kernel (bottom right).
    }
    \label{fig:results_elem1}
    \end{figure*}

\begin{figure*}[htbp]
    \centering
    \includegraphics[width=6cm]{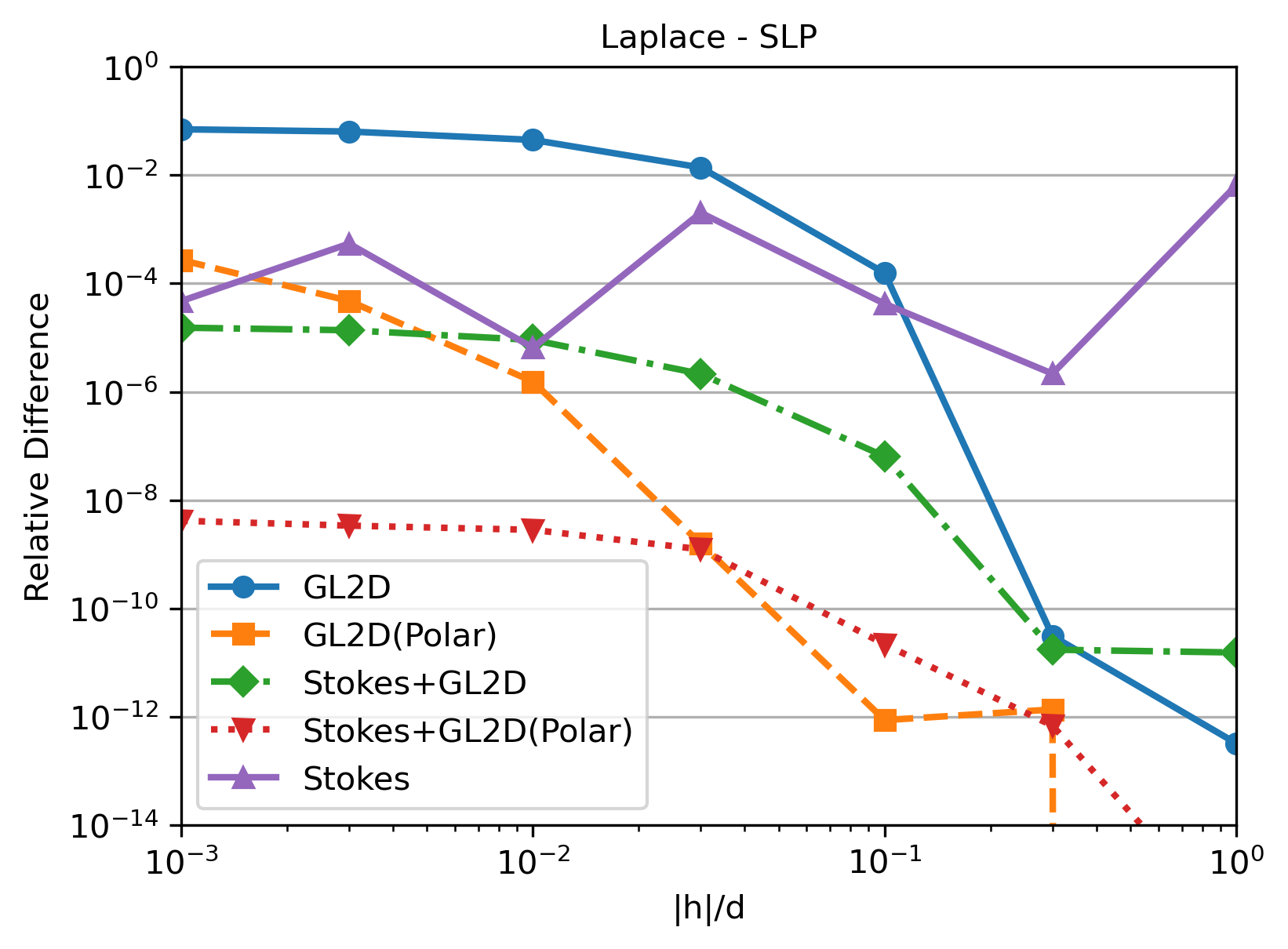}
    \includegraphics[width=6cm]{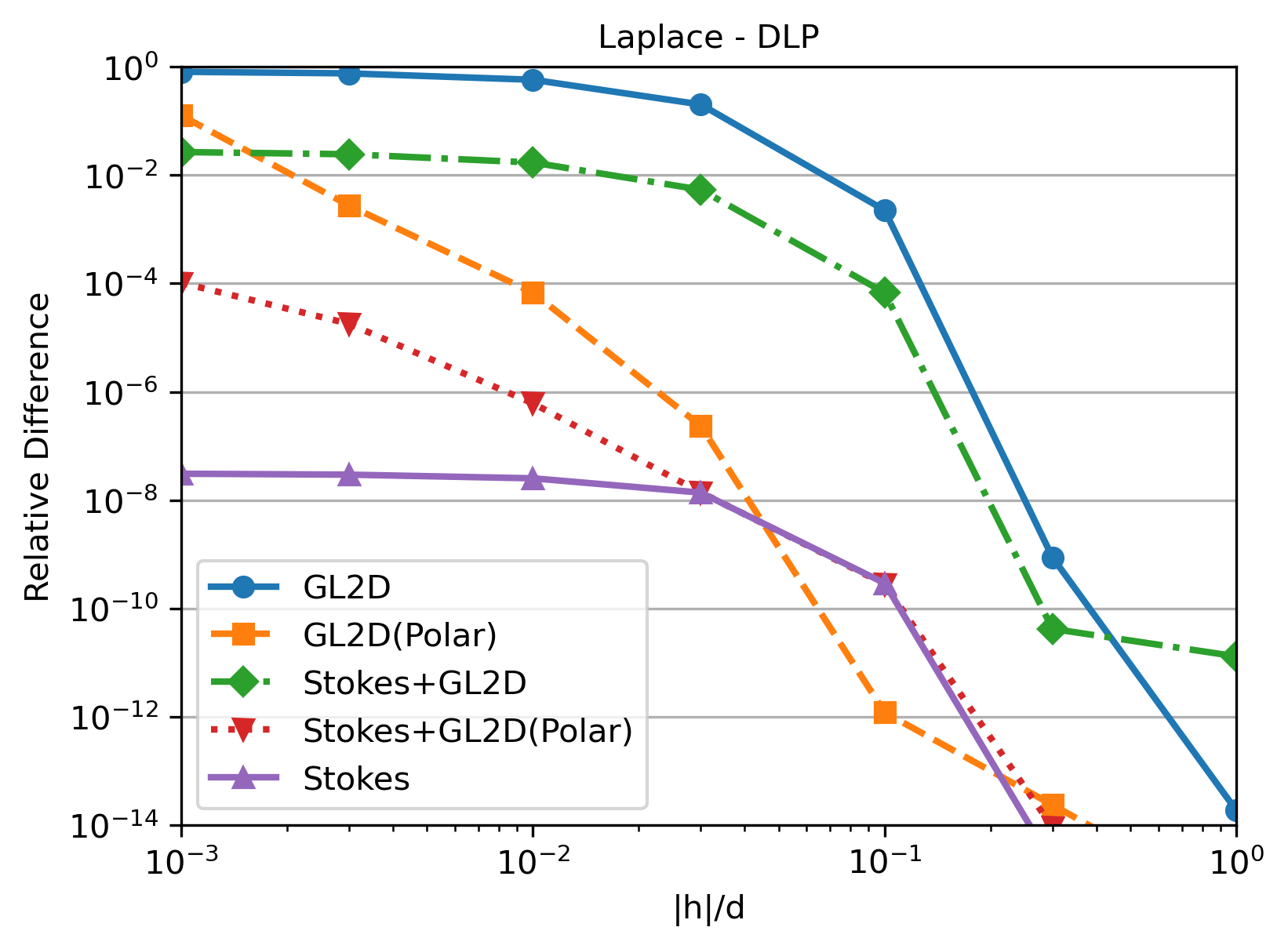}
    \includegraphics[width=6cm]{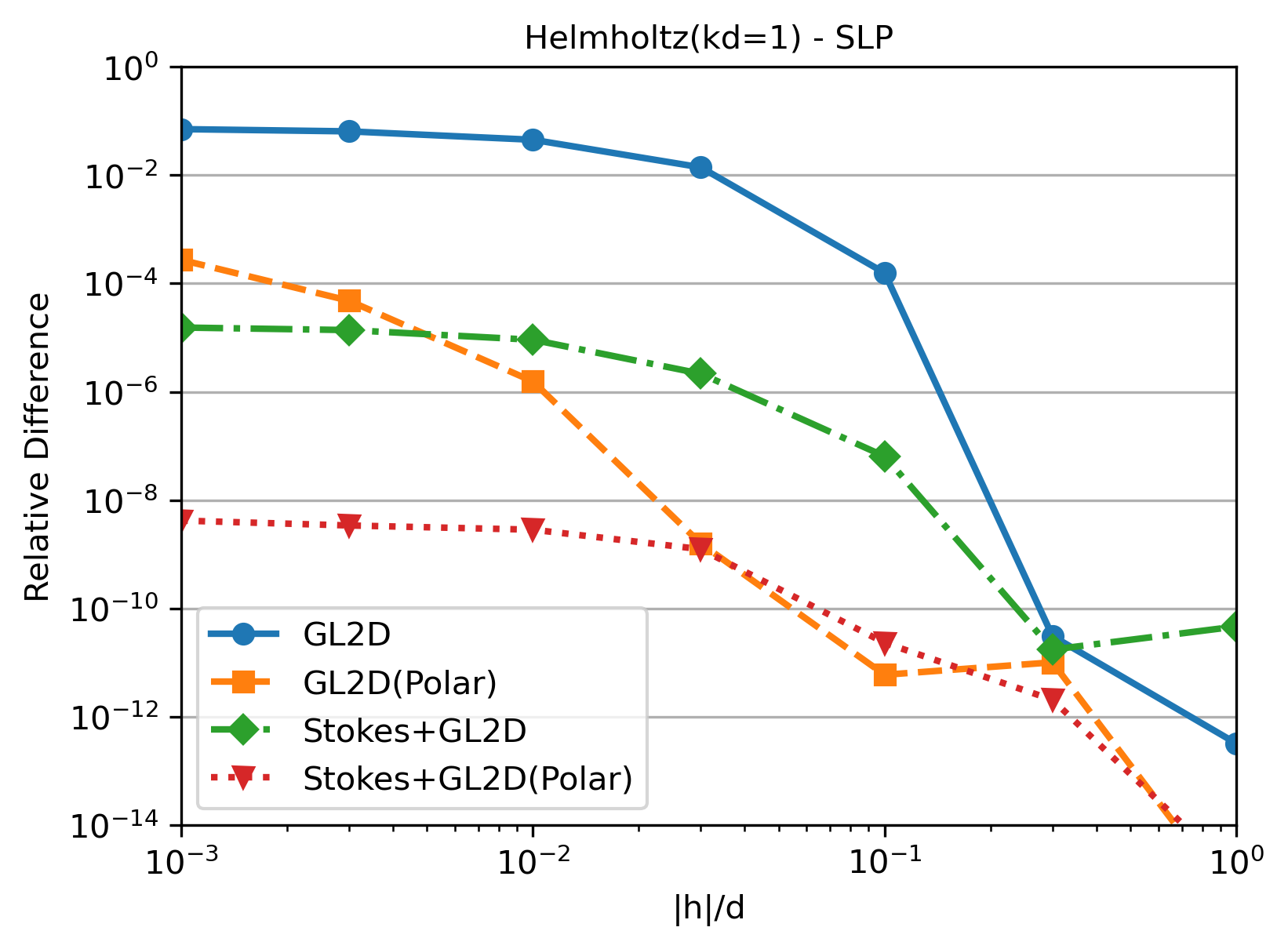}
    \includegraphics[width=6cm]{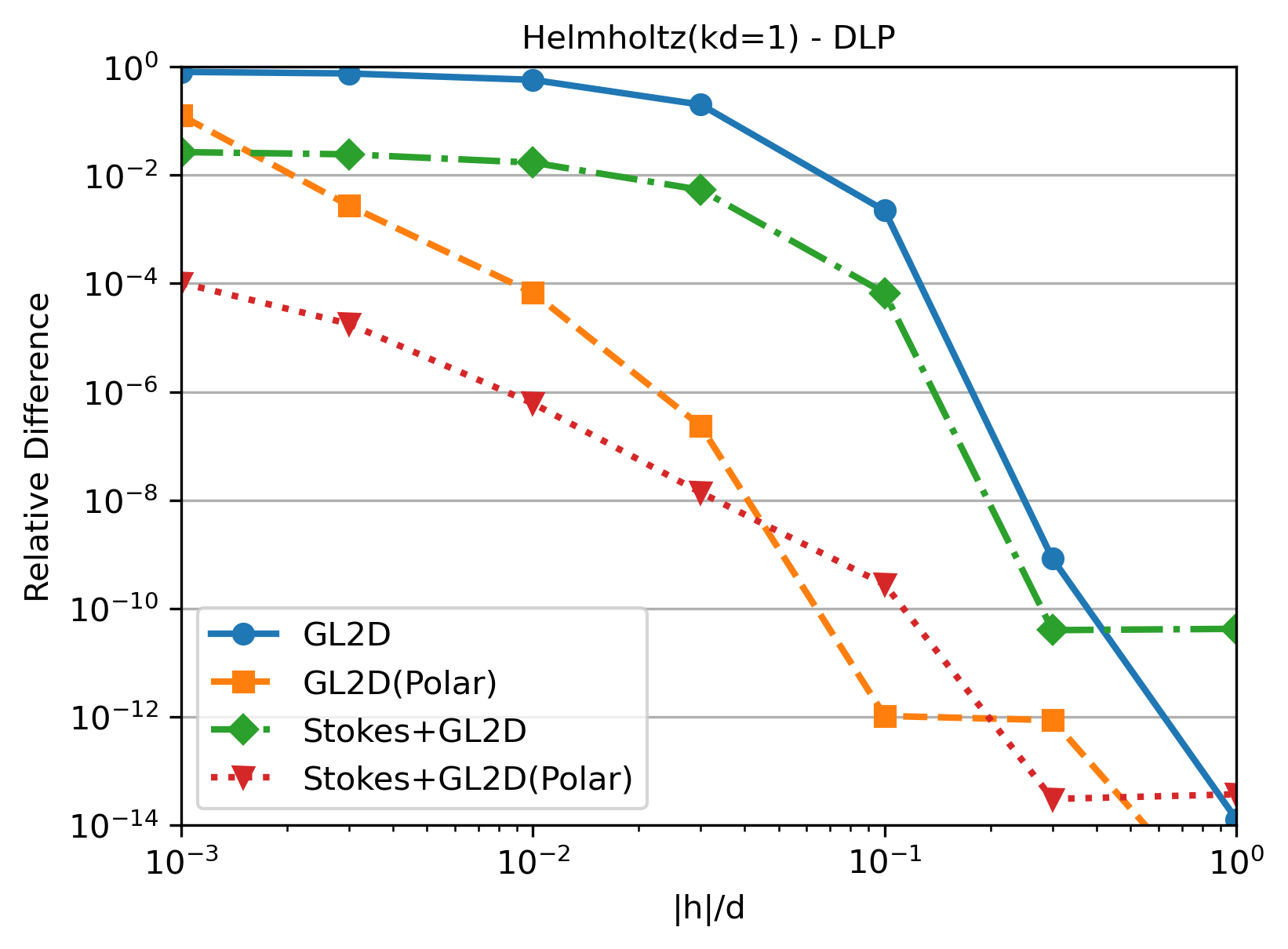}
    \caption{Relative difference of the layer potentials on \emph{element 2} for various methods and the reference adaptive Gauss-Kronrod quadrature. Results for single layer potential with Laplace kernel (top left), double layer potential with Laplace kernel (top right), single layer potential with Helmholtz kernel (bottom left), and double layer potential with Helmholtz kernel (bottom right).
    }
    \label{fig:results_elem2}
    \end{figure*}

The method was tested for both the single and double layer potentials, for both the Laplace and Helmholtz kernels. Adaptive Gauss-Kronrod quadrature, implemented in QUADPACK~\cite{piessens2012quadpack}, was used to compute the reference values of the layer potentials $P_{\mathrm{GK}}$ over a boundary element. The error tolerance of Gauss-Kronrod was set to $10^{-12}$. 
The layer potentials were computed using two standard techniques: 
(1) two-dimensional Gauss-Legendre quadrature~\cite{dunavant1985high} (\emph{GL2D}), 
(2) two-dimensional quadrature using the polar coordinate transform, i.e. nested one-dimensional Gauss-Legendre quadrature over $R$ and $\theta$ (\emph{GL2D(Polar)}), 
and with two methods using the proposed decomposition \cref{eq:greenFncDecomposition}: 
(3) evaluating the Stokes term via one-dimensional Gauss-Legendre quadrature and the curvature term via two-dimensional Gauss-Legendre quadrature~\cite{dunavant1985high} (\emph{Stokes+GL2D}), and
(4) evaluating the Stokes term via one-dimensional Gauss-Legendre quadrature and the curvature term via polar coordinate transform i.e. nested one-dimensional Gauss-Legendre quadrature  (\emph{Stokes+GL2D(Polar)}). 
For the Laplace layer potentials, we also computed the integrals using the method described in~\cite{zhu2022high} (\emph{Stokes}) using 20th degree polynomials as the basis for the density function approximation. 
The computed potentials $P$ were compared against the reference result in terms of the relative error $|P-P_{\mathrm{GK}}|/|P_{\mathrm{GK}}|$. 
20th order Gauss-Legendre quadrature was used for all integrals. 
Curved triangles parametrized via $\rq(u,v)=(u,v,f(u,v))^T$ with $f(u,v)=\sigma((u-1/4)^2+(v-1/4)^2)$ were used as test cases where $\sigma=\{-0.6,0.6\}$. 
The test elements are referred to as \emph{element 1} and \emph{element 2} and are shown in \cref{fig:testElements}.
The evaluation point was placed at $(u_p,v_p,f(u_p,v_p)+|h|\n_p)^T$, where $\n_p$ is the unit normal vector at point $\rq(u_p,v_p)$ with $(u_p,v_p)=(0.2,0.3)$. 
Results are shown in Figures \ref{fig:results_elem1} and \ref{fig:results_elem2} for \emph{elements 1} and \emph{2}, respectively. 
It was found that the proposed decomposition improves the accuracy of the numerical results in the nearly singular regime $|h|/d<1$ compared to baseline methods without the decomposition, with $d$ the maximum inter-vertex distance of the element. 
In the Laplace double layer potential case, the Stokes-only approach by~\cite{zhu2022high} (\emph{Stokes}) was found to deliver better accuracy than the proposed method. 
However, larger error was observed for the single layer potential case using this method. 
This could be because of the approximation introduced in this method where the single layer potential is considered a double layer potential with a modified density function. 
See~\cite{zhu2022high} for details on how this approximation is constructed. 
The method in~\cite{zhu2022high}, therefore, appears to be preferable for the double layer potential, while the present method may be preferable for the single layer potential as it is free from the type of approximation error introduced by the approach taken in~\cite{zhu2022high}. 
Note that the condition $kd=1$ for the Helmholtz case approximately corresponds to six wavelengths per element, which is typically used as the maximum mesh size in boundary element analysis. 
A Python implementation of the proposed method exhibited computation times comparable to the baseline \emph{GL2D(Polar)} method, although optimized implementations via compiled languages should be used for more accurate performance evaluations. 

\subsubsection{Singular case}
\begin{figure*}[htbp]
\centering
\includegraphics[width=12cm]{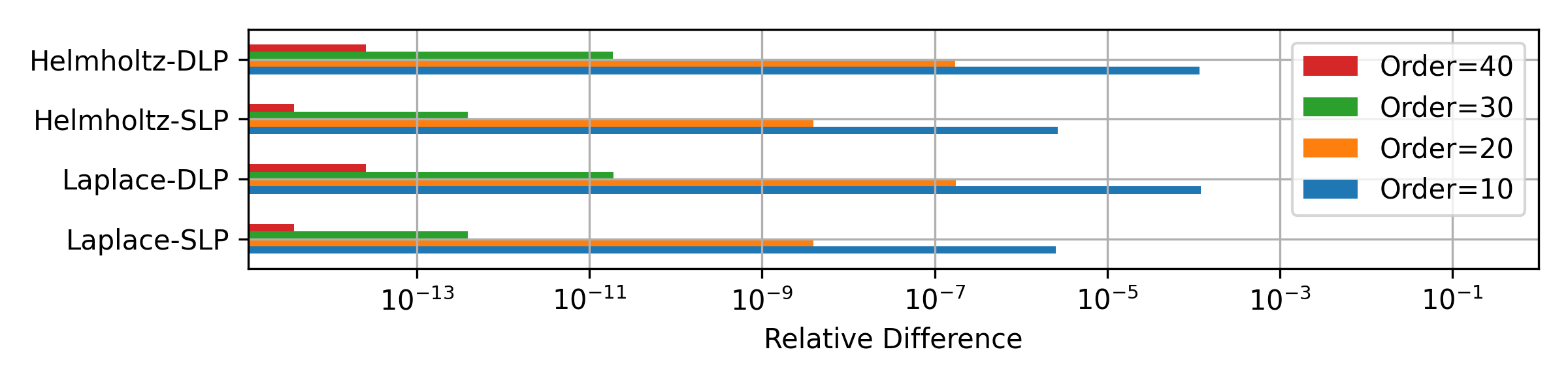}
\vspace{-15.0pt}
\caption{Relative difference $|P_{\mathrm{prop}} - P_{\mathrm{Gui}}|/|P_{\mathrm{Gui}}|$ of the singular case layer potentials on \emph{element 1} for the proposed method with various quadrature orders and the Guiggiani's reference method.}
\label{fig:rpOnElem}
\end{figure*}

Singular cases with the evaluation point $\rp$ on the element can be handled with the proposed method. 
The only change to be made is that the singularity has to be excluded from the integration domain when applying Stokes' theorem in the double layer potential case. 
Technically, this results in the subtraction of the contribution of the singularity, which is a constant of $1/2$. 
The layer potentials evaluated by the proposed method  (\emph{Stokes+GL2D(Polar)} in the previous section) $P_{\mathrm{prop}}$ were compared with reference results $P_{\mathrm{Gui}}$ obtained by Guiggiani's method~\cite{guiggiani1992general}, internally using Gauss-Legendre quadrature of 50th order. 
The quadrature order in the proposed method was varied from 10 to 40. The same quadrature order was used for the Stokes term and the curvature term. 
Element 1 from the previous section with $\rp$ on point $\rq(0.2,0.3)$ was used to compute the single and double layer potentials for the Laplace and Helmholtz kernels. 
Results in \cref{fig:rpOnElem} show $p$-convergence and a 
good agreement with the reference at sufficiently high quadrature orders. 

\subsection{Integrated BEM test: thin spherical cavity problem}

\begin{figure*}[htbp]
\centering
\includegraphics[width=5cm]{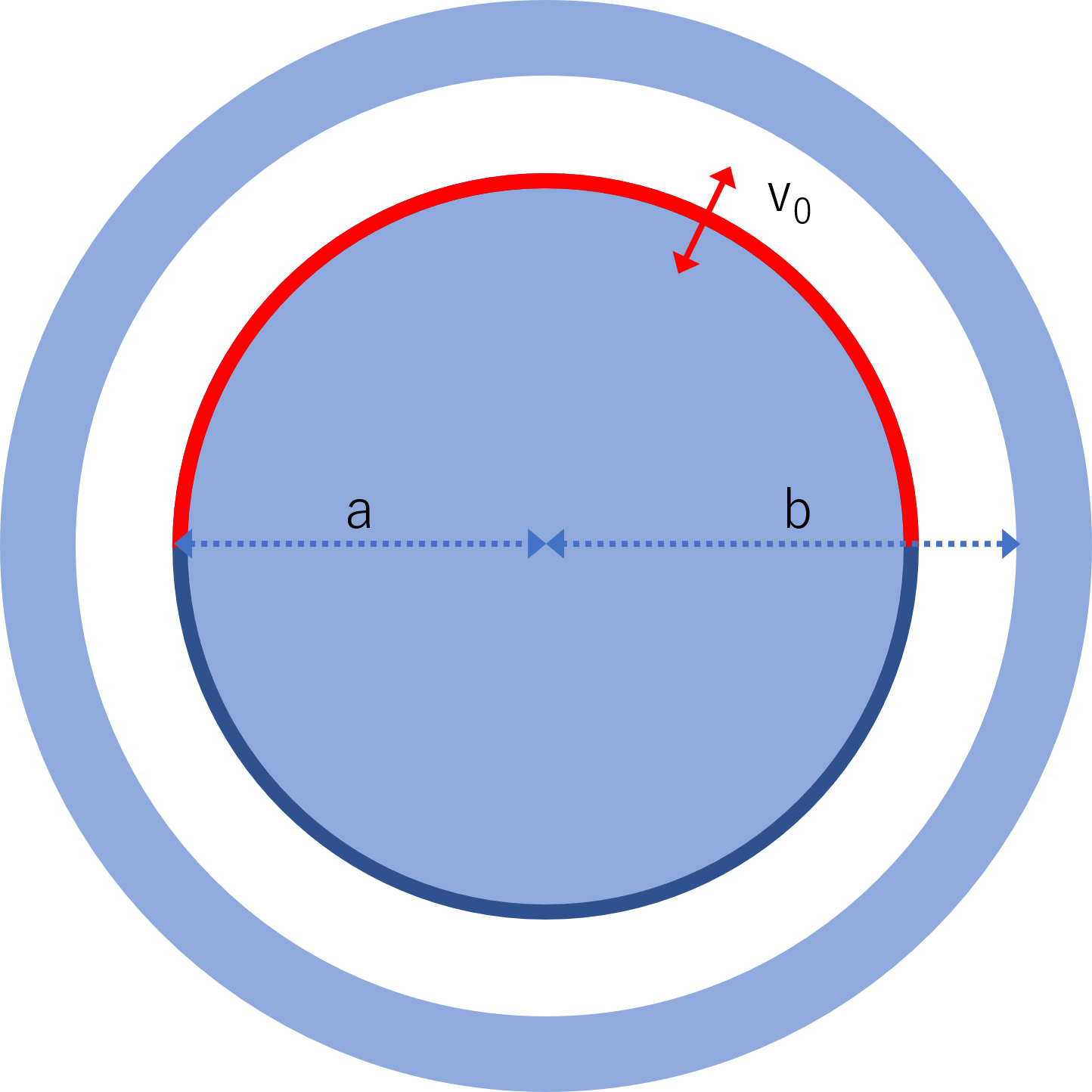}
\caption{The cross section of the spherical cavity setup used as a benchmark problem. The white region in between the two spheres is the domain $\Omega = \{r | a\le r\le b\}$ subject to analysis. The vibrating surface is the upper hemisphere of the interior sphere and is indicated as the red arc. }
\label{fig:sphericalCavity}
\end{figure*}

\begin{figure*}[htbp]
\centering
\includegraphics[width=12cm]{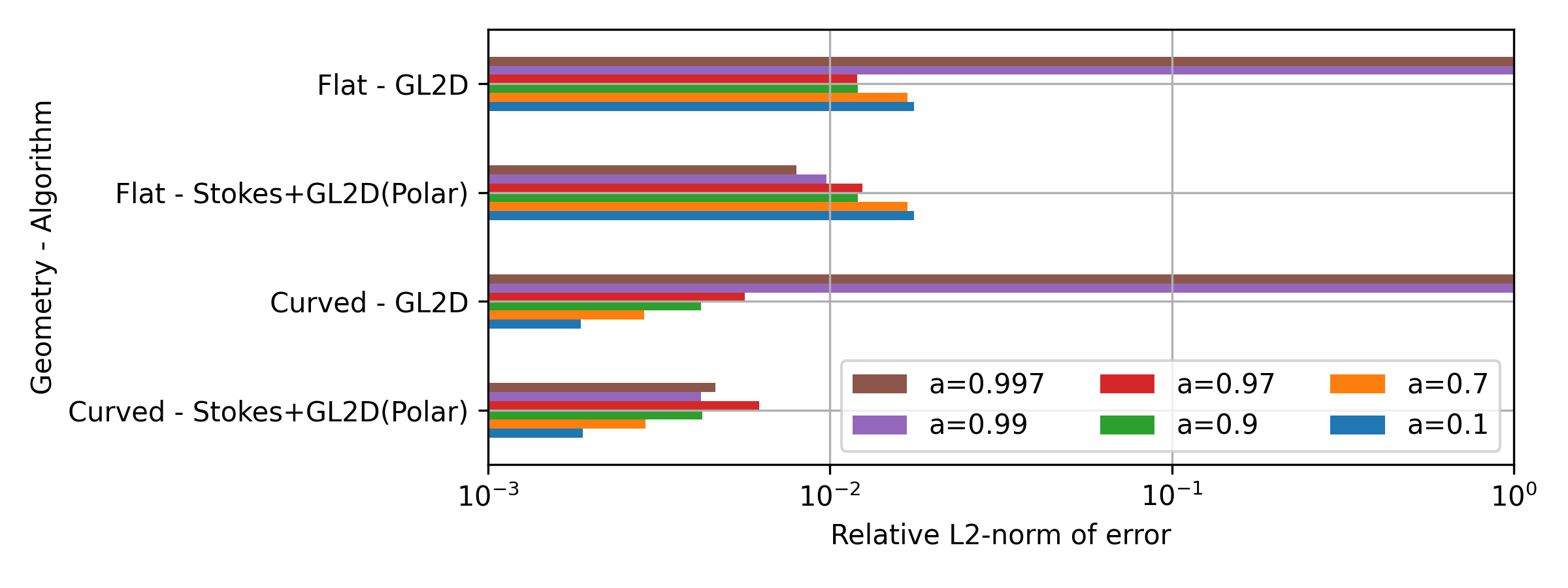}
\vspace{-15.0pt}
\caption{Results of the spherical cavity BEM test problem. Relative $L_2$-norm of the error at the collocation points for two geometry conditions (flat or non-flat manifold elements) and two algorithms used for nearly singular integrals (Gauss-Legendre quadrature or the proposed method).}
\label{fig:bemTest}
\end{figure*}

The method was integrated in a prototype BEM solver and was evaluated by solving a benchmark problem 
where we consider an interior Helmholtz problem in a spherical cavity, which is illustrated in \cref{fig:sphericalCavity}.
In this setup, a rigid sphere of radius $a$ is placed concentrically in a spherical room of radius $b$ with rigid boundaries, and the upper half $(\theta< \pi/2)$ of the internal sphere is vibrating at velocity $v_0$. 
This problem can be solved analytically and the solution expressed in spherical coordinates $(r,\theta,\phi)$ is given as follows (see e.g. \cite{williams2000fourier, gumerov2005fast}):
\begin{equation}\begin{aligned}
& p(r,\theta,\phi) = \sum_{n=0}^{\infty} \alpha_n R_n^0(r,\theta,\phi) + \beta_n S_n^0(r,\theta,\phi), \\
& R_n^m(r,\theta,\phi) = j_n(kr) Y_n^m(\theta,\phi), \quad  S_n^m(r,\theta,\phi) = h_n(kr) Y_n^m(\theta,\phi), \\
& \alpha_n = iv_0 c_{\mathrm{s}} q \sqrt{(2n+1)\pi} \left(j_n'(ka) - \frac{j_n'(kb)}{h_n'(kb)}h_n'(ka)  \right)^{-1}, \quad \beta_n = -\frac{j_n'(kb)}{h_n'(kb)}\alpha_n, \\
\end{aligned}\label{eq:analyticalSolutionOfSphericalCavity}\end{equation}
with $j_n$ and $h_n$ the spherical Bessel and Hankel function of the first kind, respectively, 
$j_n'$ and $h_n'$ the derivative of $j_n$ and $h_n$ with respect to the argument, respectively,
$Y_n^m$ the spherical harmonics, $c_{\mathrm{s}}$ the speed of sound, and $q$ the density of the medium.
We fixed $b=v_0=1$ and $k=2$. 
The numerical solution to this problem was computed using BEM where the singular integrals, i.e. layer potentials with the evaluation point on the same element, were computed using the proposed \emph{Stokes+GL2D(Polar)} method. Layer potentials with the evaluation point $\rp$ which satisfy $||\rp-\r_e|| <l_e$ with $\r_e=\rq(1/3,1/3)$ and $l_e$ the maximum length of the straight line segments connecting the vertices of the element, were considered nearly singular and were evaluated using either Gauss-Legendre quadrature (\emph{GL2D})~\cite{dunavant1985high} or the proposed \emph{Stokes+GL2D(Polar)} method. 
The quadrature order was set to 20 for all singular and nearly singular integrals. 
The manifold elements were parametrized as: 
\begin{equation}\begin{aligned}
\rq(u,v) = s \frac{\tilde{\r}(u,v)}{|\tilde{\r}(u,v)|}, \quad \tilde{\r} = \v_1 + u (\v_2-\v_1) + v (\v_3-\v_1),
\end{aligned}\end{equation}
with $s$ the radius of the spherical triangle and $\v_1, \v_2$, and $\v_3$ the vertices of the spherical triangle.
Note that this parametrization represents the surface exactly, hence the discretization error in the numerical solution for the manifold element case is solely due to the discretization of the function space. On the other hand, the flat polygon mesh introduces geometric approximation error. 
The numerical solution at the collocation points $\mathbf{p}_{\mathrm{BEM}}$ was compared against the analytical solution $\mathbf{p}_{\mathrm{exact}}$ given by \cref{eq:analyticalSolutionOfSphericalCavity} in terms of the relative $L_2$-norm of the difference vector $||\mathbf{p}_{\mathrm{BEM}} - \mathbf{p}_{\mathrm{exact}}||_2/||\mathbf{p}_{\mathrm{exact}}||_2$. 
We ran the experiments for two geometry representation conditions, where in one case the elements were represented as exact spherical triangles and in the other case polygon meshes with flat triangle elements were used to approximate the geometry of the spherical surfaces. 
In both representations, the boundary was represented by 3668 elements in total. 
The result is shown in \cref{fig:bemTest}. 
It was found that \emph{GL2D} diverges from the analytical solution in the nearly singular regime, and that the numerical results using the exact spherical surface representation delivers up to about one order of magnitude smaller error compared to the flat element counterpart. 
The error in the curved element case is bound by the constant density approximation; expanding the surface density using higher order basis functions would further improve the accuracy. 

\section{Conclusion}
A method for the evaluation of nearly singular and singular integrals for single and double layer potentials over manifold boundary elements for Laplace and Helmholtz kernels was proposed. The method uses a novel decomposition of the layer potentials into an integral of a differential 2-form which can be reduced to a contour integral via Stokes' theorem and a second term related to the curvature of the element, which can be further regularized via a polar coordinate transform and integrated via existing quadrature methods. 
Numerical tests showed that the proposed method delivers accurate results in the nearly singular and singular regime where a naïve use of Gauss-Legendre quadrature is not effective. 
One of the benefits of the proposed method, which is shared with some modern techniques (e.g.~\cite{klockner2013quadrature, morse2021robust, kaneko2023recursive}), is that it supports both nearly singular and singular layer potential integrals within the same framework. 
The proposed method covers both single and double layer potentials for both the Laplace and Helmholtz equations. 
While we introduced the method for constant densities in the present work, supporting higher order density functions is indeed desired for a boundary element method with higher order accuracy and this is a natural next step. Nevertheless, constant elements are still useful in solving large scale problems with millions or billions of unknowns under limited compute resources and in applications where accuracy requirements are relatively relaxed. 
Generalizing the present method to other layer potentials, other kernels, and integrating it into a FMM-BEM solver are other directions for future work.

\appendix
\section{Proof of decomposition \cref{eq:greenFncDecomposition}}\label{appendix_proof_decomposition}

\begin{proof}

Let us define the following \emph{pseudo potential fields} (see \cref{sec:problemStatement} and \cref{fig:elementAndCoordinateFrames} for the definition of $r$, $h$ and $\boldsymbol{\rho}$):
\begin{equation}\begin{aligned}
&\m_{\mathrm{L}} \equiv \frac{(r-h)\boldsymbol{\rho}}{4\pi \rho^2}, \quad 
\m_{\mathrm{L}}' \equiv \frac{(r-h)\boldsymbol{\rho}}{4\pi r\rho^2}, \\
&\m_{\mathrm{H}} \equiv \frac{e^{ikr}-e^{ikh}}{4\pi ik \rho^2}\boldsymbol{\rho}, \quad
\m_{\mathrm{H}}' \equiv \frac{re^{ikr}-he^{ikh}}{4\pi r\rho^2}\boldsymbol{\rho},
\end{aligned}\label{eq:mVectorFields}\end{equation}
where $\mathrm{L}$ and $\mathrm{H}$ indicate the Laplace and Helmholtz kernels respectively and the prime denotes corresponding expressions for the double layer potential. 
In the following we will use the notation:
\begin{equation}\begin{aligned}
&\r_u \equiv \frac{\partial \r_q}{\partial u}, \quad
\r_v \equiv \frac{\partial \r_q}{\partial v}, \quad
\n_u \equiv \frac{\partial \n_q}{\partial u}, \quad
\n_v \equiv \frac{\partial \n_q}{\partial v}, \\
&\r_{uu} \equiv \frac{\partial^2 \r_q}{\partial u^2}, \quad
\r_{uv} \equiv \frac{\partial^2 \r_q}{\partial u \partial v}, \quad
\r_{vv} \equiv \frac{\partial^2 \r_q}{\partial v^2}, \\
&r_{u\|} \equiv \r_u \cdot \rhohat, \quad \r_{u\|} \equiv (\r_u \cdot \rhohat) \rhohat,\quad
r_{u\perp} \equiv \r_u \cdot \rhotil, \quad \r_{u\perp} \equiv (\r_u \cdot \rhotil) \rhotil \\
&r_{v\|} \equiv \r_v \cdot \rhohat, \quad \r_{v\|} \equiv (\r_v \cdot \rhohat) \rhohat,\quad
r_{v\perp} \equiv \r_v \cdot \rhotil, \quad \r_{v\perp} \equiv (\r_v \cdot \rhotil) \rhotil\\
&n_{u\|} \equiv \n_u \cdot \rhohat, \quad \n_{u\|} \equiv (\n_u \cdot \rhohat) \rhohat,\quad
n_{u\perp} \equiv \n_u \cdot \rhotil, \quad \n_{u\perp} \equiv (\n_u \cdot \rhotil) \rhotil \\
&n_{v\|} \equiv \n_v \cdot \rhohat, \quad \n_{v\|} \equiv (\n_v \cdot \rhohat) \rhohat,\quad
n_{v\perp} \equiv \n_v \cdot \rhotil, \quad \n_{v\perp} \equiv (\n_v \cdot \rhotil) \rhotil.
\end{aligned}\end{equation}

We consider local Cartesian coordinate frames with origin at $\rq$, $x$- and $y$- axes lying in the tangent plane $T_q(S)$, and the $z$-axis pointing towards the normal direction $\n_q$ (see \cref{fig:elementAndCoordinateFrames}).
With $\boldsymbol{\rho} \equiv \rq-\rp+(\nq\cdot(\rp-\rq))\nq = \rq-\rp + h\nq$,
we have:
\begin{equation}\begin{aligned}
\partial_u \boldsymbol{\rho} &= \r_u + h \n_u + c \nq,\quad
\partial_u \rho = \boldsymbol{\hat{\rho}} \cdot \partial_u \boldsymbol{\rho} = \boldsymbol{\hat{\rho}} \cdot (\r_u + h \n_u),\quad
\partial_u r = \frac{\r_u \cdot \boldsymbol{\rho}}{r},\\
\partial_u \n_q &\equiv \n_u = \frac{1}{J} (\partial_u \mathbf{C} - \n_q (\n_q\cdot\partial_u\mathbf{C})), \quad
\mathbf{C} \equiv \r_u \times \r_v, \quad
J = |\mathbf{C}| ,\\ 
\partial_u h &= \n_u \cdot (\rp-\rq) = - \n_u \cdot \boldsymbol{\rho},\quad
\partial_z h = -1,\quad
\partial_z \boldsymbol{\rho} = \mathbf{0},\quad
\partial_z r 
= \frac{-h}{r},\\
\end{aligned}\end{equation}
where $c$ is some real number. 
It can be found that the partial derivative of the pseudo potential fields in \cref{eq:mVectorFields} with respect to $u$ and $v$ have the general form:
\begin{equation}\begin{aligned}
4\pi \partial_u \m &= A \r_{u\perp} + B \r_{u\|} + C \n_{u\perp} + D  \n_{u\|} + c \n_q, \\
4\pi \partial_v \m &= A \r_{v\perp} + B \r_{v\|} + C \n_{v\perp} + D  \n_{v\|} + c \n_q, \\
\end{aligned}\label{eq:ABCD}\end{equation}
where $\m$ is the appropriately chosen pseudo potential field in \cref{eq:mVectorFields}.
Since we have
\begin{equation}\begin{aligned}
4 \pi J \partial_x \m_x &= (\r_v)_y 4 \pi  \partial_u \m_x  -(\r_u)_y 4 \pi  \partial_v \m_x,\\
4 \pi J \partial_y \m_y &= -(\r_v)_x 4 \pi  \partial_u \m_y + (\r_u)_x 4 \pi   \partial_v \m_y, \\
\end{aligned}\end{equation}
it follows that the surface divergence $\nabla_s = \nabla - \nq(\nq\cdot\nabla)$ of $\m$ multiplied by $4\pi J$ can be computed as:
\begin{equation}\begin{aligned}\label{eq:surfaceDivergence}
& \!\!\!\! \!\!\!\! 4 \pi J (\partial_x \m_x + \partial_y \m_y) \\
  =& - (-r_{v\perp} \hat{\boldsymbol{\rho}} + r_{v\|} \tilde{\boldsymbol{\rho}}) \cdot 4 \pi  \partial_u \m  
     + (-r_{u\perp} \hat{\boldsymbol{\rho}} + r_{u\|} \tilde{\boldsymbol{\rho}}) \cdot 4 \pi    \partial_v \m\\
  =& - (-r_{v\perp} \hat{\boldsymbol{\rho}} + r_{v\|} \tilde{\boldsymbol{\rho}}) \cdot  \left( A \r_{u\perp} + B \r_{u\|} + C \n_{u\perp} + D  \n_{u\|} \right)  \\
     & + (-r_{u\perp} \hat{\boldsymbol{\rho}} + r_{u\|} \tilde{\boldsymbol{\rho}}) \cdot   \left( A \r_{v\perp} + B \r_{v\|} + C \n_{v\perp} + D  \n_{v\|} \right) \\
 =& +r_{v\perp} \left(  B r_{u\|} + D   n_{u\|} \right) - r_{v\|}  \left( A r_{u\perp} + C n_{u\perp} \right) \\
  & - r_{u\perp} \left( B r_{v\|} + D n_{v\|} \right) + r_{u\|} \left( A r_{v\perp} + C n_{v\perp} \right) \\
 =& (A+B)(r_{u\|} r_{v\perp} - r_{u\perp} r_{v\|} ) + C ( r_{u\|}    n_{v\perp} - r_{v\|}     n_{u\perp} ) 
      + D ( r_{v\perp}   n_{u\|}    - r_{u\perp}  n_{v\|}   ) \\
 =& (A+B)J - C J \kappa_N(\rhotil) - D J \kappa_N(\rhohat), \\
\end{aligned}\end{equation}
where $\kappa_N(\rhotil)$ and $\kappa_N(\rhohat)$ are the normal curvatures at point $\rq$ on the element with respect to normal planes spanned by $\n_q$ and tangent vectors $\rhotil$ and $\rhohat$, respectively. 
In the last step, we have used the following lemma:
\begin{lemma}
\begin{equation}\begin{aligned}
   C ( r_{v\|}     n_{u\perp} - r_{u\|}     n_{v\perp}  ) 
 + D ( r_{u\perp} n_{v\|}    - r_{v\perp} n_{u\|}     ) 
 =  J(C \kappa_N(\rhotil) + D  \kappa_N(\rhohat)), \\
\end{aligned}\label{eq:relationOf_ru_rv_nu_nv_and_kappa}\end{equation}
\end{lemma}
\begin{proof}
See \cref{appendix_normal_curvature}.
\end{proof}
From \cref{eq:surfaceDivergence} it follows:
\begin{equation}\begin{aligned}
\frac{A+B}{4\pi} =  \nabla_s \cdot \m  + \frac{C\kappa_N(\rhotil)+D\kappa_N(\rhohat)}{4\pi}. \\
\end{aligned}\label{eq:oneStepBeforeDecomposition}\end{equation}
The coefficients $A+B$, $C$, and $D$ for each pseudo potential field in \cref{eq:mVectorFields} are summarized in \cref{table:expressionsForQBC}.
\begin{table*}
\caption{Coefficients introduced in \cref{eq:ABCD} for Laplace and Helmholtz, single and double layer potentials.}
\footnotesize
\begin{equation}\label{table:expressionsForQBC}
\begin{array}{c||c| c | c | c |c|c} 
& A+B
& C 
& D 
\\\hline 
\m_{\mathrm{L}}
& \frac{1}{r}
& \frac{h}{r+h}
& \frac{r}{r+h}
\\\hline 
\m_{\mathrm{L}}' 
& \frac{h}{r^3}
& \frac{h}{r(r+h)}
& \frac{1}{r+h}
\\\hline 
\m_{\mathrm{H}} 
& \frac{e^{ikr}}{r}
& \frac{h(e^{ikr}-e^{ikh})}{ik\rho^2}
& e^{ikh} - \frac{h(e^{ikr}-e^{ikh})}{ik\rho^2}
\\\hline 
\m_{\mathrm{H}}' 
& \frac{he^{ikr}(1-ikr)}{r^3}
& \frac{h\left(r e^{i k h}-h e^{i k r}\right)}{r \rho^2} 
& \frac{r e^{i k r}-h e^{i k h}}{ \rho^2}-i k e^{i k h}
\\
\end{array}
\end{equation}
\vspace{-15.0pt}
\end{table*}
It turns out that for all cases, the left hand side of \cref{eq:oneStepBeforeDecomposition} is nothing but the Green function or its normal derivative. 
Lastly, by using $\mathbf{f}_\mathrm{K} = \n_q \times \m_{\mathrm{K}}$ and $\mathbf{f}_\mathrm{K}' = \n_q \times \m_{\mathrm{K}}'$, we finally obtain decomposition \cref{eq:greenFncDecomposition}:
\begin{equation}\begin{aligned}
G_{\mathrm{K}}(\rp,\rq) &=  (\nabla_{\rq} \times  \mathbf{f}_{\mathrm{K}} )\cdot \nq  + \frac{1}{4\pi}\left(C_{\mathrm{K}}\kappa_N(\rhotil)+D_{\mathrm{K}}\kappa_N(\rhohat)\right),\\
\frac{\partial G_{\mathrm{K}}(\rp,\rq)}{\partial \nq} &=  (\nabla_{\rq} \times  \mathbf{f}_{\mathrm{K}}' )\cdot \nq + \frac{1}{4\pi}\left(C'_{\mathrm{K}}\kappa_N(\rhotil)+D'_{\mathrm{K}}\kappa_N(\rhohat)\right).\\
\end{aligned}\end{equation}

\end{proof}

\section{Proof of relation \cref{eq:relationOf_ru_rv_nu_nv_and_kappa}}\label{appendix_normal_curvature}

We use the following definitions of first and second fundamental forms~\cite{do2016differential}:
\begin{equation}\begin{aligned}
E \equiv \r_u^2, \ \ F \equiv \r_u\cdot\r_v, \ \ G \equiv \r_v^2, \ \
e \equiv \r_{uu}\cdot \nq, \ \ f \equiv \r_{uv}\cdot\nq, \ \ g \equiv \r_{vv}\cdot\nq.\\
\end{aligned}\end{equation}
With $\mathbf{C}_u \equiv \partial \mathbf{C}/\partial u$,  $\mathbf{C}_v \equiv \partial \mathbf{C}/\partial v$, $\theta_u$ the angle of $\rhohat$ from $\r_u$ and $\theta_v$ the angle of $\rhohat$ from $\r_v$ measured in the tangent plane of $S$ at $\rq$, 
\begin{equation}\begin{aligned}
K \equiv& - D r_{u\perp}  n_{v\|}
      + Cr_{u\|}     n_{v\perp}
      + Dr_{v\perp}  n_{u\|}
      - Cr_{v\|}     n_{u\perp} \\
 =&   - D(\r_{u}\cdot\rhotil)  (\n_{v}\cdot\rhohat) 
      + C(\r_{u}\cdot\rhohat)  (\n_{v}\cdot\rhotil)\\
    &  + D(\r_{v}\cdot\rhotil)  (\n_{u}\cdot\rhohat)  
      - C(\r_{v}\cdot\rhohat)  (\n_{u}\cdot\rhotil) \\
 =&   - D(\r_{u}\cdot\rhotil)  (\n_{v}\cdot\rhohat) 
      + D(\r_{u}\cdot\rhohat)  (\n_{v}\cdot\rhotil)\\
    &  + C(\r_{v}\cdot\rhotil)  (\n_{u}\cdot\rhohat)  
      - C(\r_{v}\cdot\rhohat)  (\n_{u}\cdot\rhotil) \\
   &+  (C-D)(\r_{u}\cdot\rhohat)  (\n_{v}\cdot\rhotil)
    -  (C-D)(\r_{v}\cdot\rhotil)  (\n_{u}\cdot\rhohat) \\
 =&  C(\r_v \times \n_u) \cdot (\rhotil \times \rhohat)
    - D(\r_u \times \n_v) \cdot (\rhotil \times \rhohat) \\
    &+ (C-D) \left( (\r_{u}\cdot\rhohat)  (\n_{v}\cdot\rhotil)
                 - (\r_{v}\cdot\rhotil)  (\n_{u}\cdot\rhohat) \right) \\
 =&  -  \nq \cdot (-D(\r_u \times \n_v) + C(\r_v \times \n_u))  \\
  & + (C-D) \left( (\r_{u}\cdot\rhohat)  (\n_{v}\cdot\rhotil)
    - (\r_{v}\cdot\rhotil)  (\n_{u}\cdot\rhohat) \right) \\
 =&    \frac{1}{J} \nq \cdot (D(\r_u \times \mathbf{C}_v) - C(\r_v \times \mathbf{C}_u))  \\
   &+ (C-D) \left( (\r_{u}\cdot\rhohat)  (\n_{v}\cdot\rhotil)
    - (\r_{v}\cdot\rhotil)  (\n_{u}\cdot\rhohat) \right) \\
 =&  -  \frac{1}{J}(C |\r_v|^2 (\r_{uu}\cdot \nq) + D |\r_u|^2 ( \r_{vv}\cdot \nq)  - (C+D)(\r_u\cdot\r_v)(\r_{uv}\cdot\nq)  )\\
    &+ (C-D) \left( (\r_{u}\cdot\rhohat)  (\n_{v}\cdot\rhotil)
    - (\r_{v}\cdot\rhotil)  (\n_{u}\cdot\rhohat) \right) \\
 =& -  \frac{1}{J}(C|\r_v|^2 (\r_{uu}\cdot \nq) + D|\r_u|^2 ( \r_{vv}\cdot \nq)  - (C+D)(\r_u\cdot\r_v)(\r_{uv}\cdot\nq)   )  \\
  & + \frac{C-D}{J} |\r_{u}|\cos\theta_u  (|\r_v|\cos\theta_v \r_{uv}\cdot\nq -|\r_u|\cos\theta_u \r_{vv}\cdot\n)\\
  & - \frac{C-D}{J} |\r_{v}|\sin\theta_v  (-|\r_v|\sin\theta_v \r_{uu}\cdot\nq + |\r_u|\sin\theta_u \r_{uv}\cdot\nq). \\
\end{aligned}\label{eq:long_1}\end{equation}
Here we used:
\begin{equation}\begin{aligned}
&\r_u\cdot\rhohat=|\r_u|\cos\theta_u, \quad \r_u\cdot\rhotil=|\r_u|\sin\theta_u, \\
&\r_v\cdot\rhohat=|\r_v|\cos\theta_v, \quad \r_v\cdot\rhotil=|\r_v|\sin\theta_v. 
\end{aligned}\end{equation}
\cref{eq:long_1} continues as:
\begin{equation}\begin{aligned}
K
 =& -  \frac{1}{J}(C|\r_v|^2 (\r_{uu}\cdot \nq) + D|\r_u|^2 ( \r_{vv}\cdot \nq)  - (C+D)(\r_u\cdot\r_v)(\r_{uv}\cdot\nq)   )  \\
  & +  \frac{C-D}{J}\left( 
     \sin^2\theta_v  |\r_{v}|^2 (\r_{uu}\cdot\nq)
    - \cos^2\theta_u |\r_{u}|^2 (\r_{vv}\cdot\nq)
    \right) \\
  & +  \frac{C-D}{J}|\r_u||\r_v| \cos(\theta_u+\theta_v) (\r_{uv}\cdot\nq) \\
 =&   \frac{1}{J}\left( 
     ((C-D)\sin^2\theta_v-C)  |\r_{v}|^2 (\r_{uu}\cdot\nq) \right) \\
   & +  \frac{1}{J}\left( (-D - (C-D)\cos^2\theta_u) |\r_{u}|^2 (\r_{vv}\cdot\nq)
    \right) \\
  & + 2\frac{1}{J}(C\cos\theta_u\cos\theta_v-D\sin\theta_u\sin\theta_v)  |\r_u||\r_v| (\r_{uv}\cdot\nq) \\
 =&   \frac{-1}{J}\left( 
     (D\sin^2\theta_v + C\cos^2\theta_v)  |\r_{v}|^2 (\r_{uu}\cdot\nq)\right)\\
   & + \frac{-1}{J}\left( (D\sin^2\theta_u + C\cos^2\theta_u) |\r_{u}|^2 (\r_{vv}\cdot\nq)
    \right) \\
  & + 2\frac{1}{J}(C\cos\theta_u\cos\theta_v-D\sin\theta_u\sin\theta_v)  |\r_u||\r_v| (\r_{uv}\cdot\nq) \\
 =&  - \frac{D}{J}\left( 
     \sin^2\theta_v  \r_{v}^2 e
    + \sin^2\theta_u \r_{u}^2 g
     - 2\sin\theta_u\sin\theta_v  |\r_u||\r_v| f \right) \\
 & -  \frac{C}{J}\left( 
     \cos^2\theta_v  \r_{v}^2 e
    + \cos^2\theta_u \r_{u}^2 g
   - 2\cos\theta_u\cos\theta_v |\r_u||\r_v| f \right)   \\
 =&   -JD \frac{
      (\sin^2\theta_v) G e
    + (\sin^2\theta_u) E g
    - 2(\sin\theta_u\sin\theta_v)  |\r_u||\r_v| f }{EG-F^2} \\
 & -  JC  \frac{
      (\cos^2\theta_v) G e
    + (\cos^2\theta_u) E g
    - 2(\cos\theta_u\cos\theta_v)  |\r_u||\r_v| f }{EG-F^2}  \\
 =&   - J (C \kappa_{N}(\rhotil) + D \kappa_{N}(\rhohat) ).\\
\end{aligned}\end{equation}

\begin{remark}
From the definition of the normal curvature, it follows that the quantities $\kappa_N(\rhohat)$ and $\kappa_N(\rhotil)$ given by:
\begin{equation}\begin{aligned}
\kappa_N(\rhohat) 
     =& \frac{Ge\sin^2\theta_v + Eg\sin^2\theta_u - 2f |\r_u||\r_v| \sin\theta_u\sin\theta_v }{EG-F^2},\\
\kappa_N(\rhotil)
     =& \frac{Ge\cos^2\theta_v + Eg\cos^2\theta_u - 2f |\r_u||\r_v| \cos\theta_u\cos\theta_v }{EG-F^2},\\
\end{aligned}\end{equation}
are nothing but the normal curvature of the surface at $\rq$ in direction $\rhohat$ and $\rhotil$, respectively.
\end{remark}

\section{Acknowledgments}
This work is supported by Cooperative Research Agreement W911NF2020213 between the University of Maryland and the Army Research Laboratory, with David Hull and Steven Vinci as Technical monitors. 
Shoken Kaneko acknowledges scholarships from Japan Student Services Organization and Watanabe Foundation.
The authors would like to thank the anonymous reviewers who helped improving the manuscript.

\bibliographystyle{plain}
\bibliography{references}

\begin{thebibliography}{10}

\bibitem{adelman2016computation}
Ross Adelman, Nail~A Gumerov, and Ramani Duraiswami.
\newblock Computation of {G}alerkin double surface integrals in the 3-{D}
  boundary element method.
\newblock {\em IEEE Trans. Antennas Propag.}, 64(6):2389--2400, 2016.

\bibitem{beer2020isogeometric}
Gernot Beer, Benjamin Marussig, and Christian Duenser.
\newblock {\em The isogeometric boundary element method}.
\newblock Springer, 2020.

\bibitem{do2016differential}
Manfredo~P Do~Carmo.
\newblock {\em Differential geometry of curves and surfaces: revised and
  updated second edition}.
\newblock Courier Dover Publications, 2016.

\bibitem{dunavant1985high}
D.~Dunavant.
\newblock High degree efficient symmetrical gaussian quadrature rules for the
  triangle.
\newblock {\em International journal for numerical methods in engineering},
  21(6):1129--1148, 1985.

\bibitem{greengard2021fast}
Leslie Greengard, Michael O'Neil, Manas Rachh, and Felipe Vico.
\newblock Fast multipole methods for the evaluation of layer potentials with
  locally-corrected quadratures.
\newblock {\em J. Comput. Phys.: X}, 10:100092, 2021.

\bibitem{guiggiani1992general}
M.~Guiggiani, G.~Krishnasamy, T.~J. Rudolphi, and F.~J. Rizzo.
\newblock A {G}eneral {A}lgorithm for the {N}umerical {S}olution of
  {H}ypersingular {B}oundary {I}ntegral {E}quations.
\newblock {\em J. Appl. Mech.}, 59(3):604--614, 09 1992.

\bibitem{gumerov2005fast}
Nail~A Gumerov and Ramani Duraiswami.
\newblock {\em Fast multipole methods for the {H}elmholtz equation in three
  dimensions}.
\newblock Elsevier, 2005.

\bibitem{gumerov2021analytical}
Nail~A Gumerov and Ramani Duraiswami.
\newblock Analytical computation of boundary integrals for the {H}elmholtz
  equation in three dimensions.
\newblock {\em arXiv:2103.17196}, 2021.

\bibitem{gumerov2023analytical}
Nail~A Gumerov, Shoken Kaneko, and Ramani Duraiswami.
\newblock Analytical {G}alerkin boundary integrals of {L}aplace kernel layer
  potentials in $\mathbb{R}^{3}$.
\newblock {\em arXiv preprint arXiv:2302.03247}, 2023.

\bibitem{GUMEROV2023recursive}
Nail~A Gumerov, Shoken Kaneko, and Ramani Duraiswami.
\newblock Recursive {C}omputation of the {M}ultipole {E}xpansions of {L}ayer
  {P}otential {I}ntegrals over {S}implices for {E}fficient {F}ast {M}ultipole
  {A}ccelerated {B}oundary {E}lements.
\newblock {\em Journal of Computational Physics}, 486(1):112118, 2023.

\bibitem{hackbusch1994numerical}
Wolfgang Hackbusch and Stefan~A Sauter.
\newblock On numerical cubatures of nearly singular surface integrals arising
  in {BEM} collocation.
\newblock {\em Computing}, 52(2):139--159, 1994.

\bibitem{hayami2005variable}
Ken Hayami.
\newblock Variable transformations for nearly singular integrals in the
  boundary element method.
\newblock {\em Publications of the Research Institute for Mathematical
  Sciences}, 41(4):821--842, 2005.

\bibitem{hayami1988quadrature}
Ken Hayami and CA~Brebbia.
\newblock Quadrature methods for singular and nearly singular integrals in
  3-{D} boundary element method.
\newblock {\em Boundary elements X}, 1:237--264, 1988.

\bibitem{johnston2007sinh}
Barbara~M Johnston, Peter~R Johnston, and David Elliott.
\newblock A sinh transformation for evaluating two-dimensional nearly singular
  boundary element integrals.
\newblock {\em Int. J. Numer. Methods Eng.}, 69(7):1460--1479, 2007.

\bibitem{kaneko2023efficient}
Shoken Kaneko and Ramani Duraiswami.
\newblock Efficient {E}xact {Q}uadrature of {R}egular {S}olid {H}armonics
  {T}imes {P}olynomials {O}ver {S}implices in $\mathbb{R}^3$.
\newblock {\em arXiv preprint arXiv:2307.12202}, 2023.

\bibitem{kaneko2023recursive}
Shoken Kaneko, Nail~A Gumerov, and Ramani Duraiswami.
\newblock Recursive {A}nalytical {Q}uadrature of {L}aplace and {H}elmholtz
  {L}ayer {P}otentials in $\mathbb{R}^{3}$.
\newblock {\em arXiv preprint arXiv:2302.02196}, 2023.

\bibitem{klockner2013quadrature}
Andreas Kl{\"o}ckner, Alexander Barnett, Leslie Greengard, and Michael O'Neil.
\newblock Quadrature by expansion: {A} new method for the evaluation of layer
  potentials.
\newblock {\em J. Comput. Phys.}, 252:332--349, 2013.

\bibitem{lee2018introduction}
John~M Lee.
\newblock {\em Introduction to {R}iemannian manifolds}, volume~2.
\newblock Springer, 2018.

\bibitem{lenoir2012evaluation}
Marc Lenoir and Nicolas Salles.
\newblock Evaluation of 3-{D} singular and nearly singular integrals in
  {G}alerkin {BEM} for thin layers.
\newblock {\em SIAM J. Sci. Comput.}, 34(6):A3057--A3078, 2012.

\bibitem{montanelli2022computing}
Hadrien Montanelli, Matthieu Aussal, and Houssem Haddar.
\newblock Computing weakly singular and near-singular integrals over curved
  boundary elements.
\newblock {\em SIAM Journal on Scientific Computing}, 44(6):A3728--A3753, 2022.

\bibitem{morse2021robust}
Matthew~J Morse, Abtin Rahimian, and Denis Zorin.
\newblock A robust solver for elliptic {PDE}s in 3{D} complex geometries.
\newblock {\em Journal of Computational Physics}, 442:110511, 2021.

\bibitem{newman1986distributions}
John~Nicholas Newman.
\newblock Distributions of sources and normal dipoles over a quadrilateral
  panel.
\newblock {\em Journal of Engineering Mathematics}, 20(2):113--126, 1986.

\bibitem{piessens2012quadpack}
Robert Piessens, Elise de~Doncker-Kapenga, Christoph~W {\"U}berhuber, and
  David~K Kahaner.
\newblock {\em Quadpack: a subroutine package for automatic integration},
  volume~1.
\newblock Springer, 2012.

\bibitem{rosen1995continuation}
Dan Rosen and Donald~E Cormack.
\newblock The continuation approach: {A} general framework for the analysis and
  evaluation of singular and near-singular integrals.
\newblock {\em SIAM Journal on Applied Mathematics}, 55(3):723--762, 1995.

\bibitem{sauter2010boundary}
Stefan~A Sauter and Christoph Schwab.
\newblock Boundary element methods.
\newblock In {\em Boundary {E}lement {M}ethods}, pages 183--287. Springer,
  2010.

\bibitem{spivak2018calculus}
Michael Spivak.
\newblock {\em Calculus on manifolds: a modern approach to classical theorems
  of advanced calculus}.
\newblock CRC press, 2018.

\bibitem{wala2019fast}
Matt Wala and Andreas Kl{\"o}ckner.
\newblock A fast algorithm for quadrature by expansion in three dimensions.
\newblock {\em J. Comput. Phys.}, 388:655--689, 2019.

\bibitem{wala2020optimization}
Matt Wala and Andreas Kl{\"o}ckner.
\newblock Optimization of fast algorithms for global {Q}uadrature by
  {E}xpansion using target-specific expansions.
\newblock {\em J. Comput. Phys.}, 403:108976, 2020.

\bibitem{williams2000fourier}
Earl~G Williams and J~Adin Mann~III.
\newblock Fourier acoustics: sound radiation and nearfield acoustical
  holography, 2000.

\bibitem{zhu2022high}
Hai Zhu and Shravan Veerapaneni.
\newblock High-order close evaluation of {L}aplace layer potentials: {A}
  differential geometric approach.
\newblock {\em SIAM J. Sci. Comput.}, 44(3):A1381--A1404, 2022.

\end{thebibliography}
\end{document}